\RequirePackage{fix-cm}
\RequirePackage{snapshot} 

\documentclass[twocolumn]{svjour3}
\def\makeheadbox{}

\smartqed 
\usepackage{graphicx}

\journalname{}

\newcommand{\m}[1]{\mathbf{#1}} 
\newcommand{\vt}[1]{\mathbf{#1}} 
\newcommand{\vg}[1]{\bm{#1}} 
\newcommand{\T}[0]{\mathrm{T}}

\newcommand{\Ez}[0]{\mathbf{E}_0}
\newcommand{\Eo}[0]{\mathbf{E}_1}
\newcommand{\Et}[0]{\mathbf{E}_2}

\newcommand{\M}[0]{\mathbf{M}_0}
\newcommand{\K}[0]{\mathbf{K}}

\usepackage[Symbol]{upgreek}
\usepackage{bm}
\usepackage{epstopdf}
\usepackage{nicefrac}
\usepackage{lineno}
\usepackage[hidelinks]{hyperref}
\modulolinenumbers[5]
\usepackage{amsmath,amssymb}
\usepackage{subfig}
\usepackage{placeins}
\usepackage{color}

\begin{document}

\title{Automatic 3D modeling by combining SBFEM and transfinite element shape functions}
\subtitle{}

\titlerunning{Automatic 3D modeling by combining SBFEM and transfinite element shape functions}

\author{Hauke Gravenkamp  \and
    Albert A.~Saputra \and
    Sascha Eisenträger
}
\institute{H.~Gravenkamp (corresponding author)\at
    \email{hgravenkamp@cimne.upc.edu}
    }

\date{preprint submitted: 22 Oct 2019\\[2mm] 
final version published as: "Three-dimensional image-based modeling by combining SBFEM and transfinite element shape functions", Computational Mechanics (2020) 66:911–930, 10.1007/s00466-020-01884-4}

\maketitle

\begin{abstract}
    The scaled boundary finite element method (SBFEM) has recently been employed as an efficient means to model three-dimensional structures, in particular when the geometry is provided as a voxel-based image. To this end, an octree decomposition of the computational domain is deployed and each cubic cell is treated as an SBFEM subdomain. The surfaces of each subdomain are discretized in the finite element sense. We improve on this idea by combining the semi-analytical concept of the SBFEM with certain transition elements on the subdomains' surfaces. Thus, we avoid the triangulation of surfaces employed in previous works and consequently reduce the number of surface elements and degrees of freedom. In addition, these discretizations allow coupling elements of arbitrary order such that local $p$-refinement can be achieved straightforwardly.
\end{abstract}

\section{Introduction}
\label{sec:Intro}
The scaled boundary finite element method (SBFEM) is a semi-analytical technique -- loosely based on finite elements -- that involves a discretization of boundaries of computational (sub-)domains. Roughly speaking, this method aims at transforming a partial differential equation (PDE) in two or three spatial coordinates into a set of ordinary differential equations (ODE) in one coordinate by discretizing all but this one coordinate. In order to apply this idea effectively, a particular coordinate system is usually chosen in which one coordinate $\xi$ points from the origin\footnote{In the context of the SBFEM, the origin of the coordinate system is usually positioned inside the domain and referred to as \textit{scaling center}.} to the boundary while the remaining one or two coordinate(s) ($\eta$, $\zeta$) describe a parametrization of the boundary.

\begin{figure}[b!]
    \centering
    \subfloat[Quadtree mesh, SBFEM \label{fig:sketchQuad}]{\includegraphics[height=0.23\textwidth]{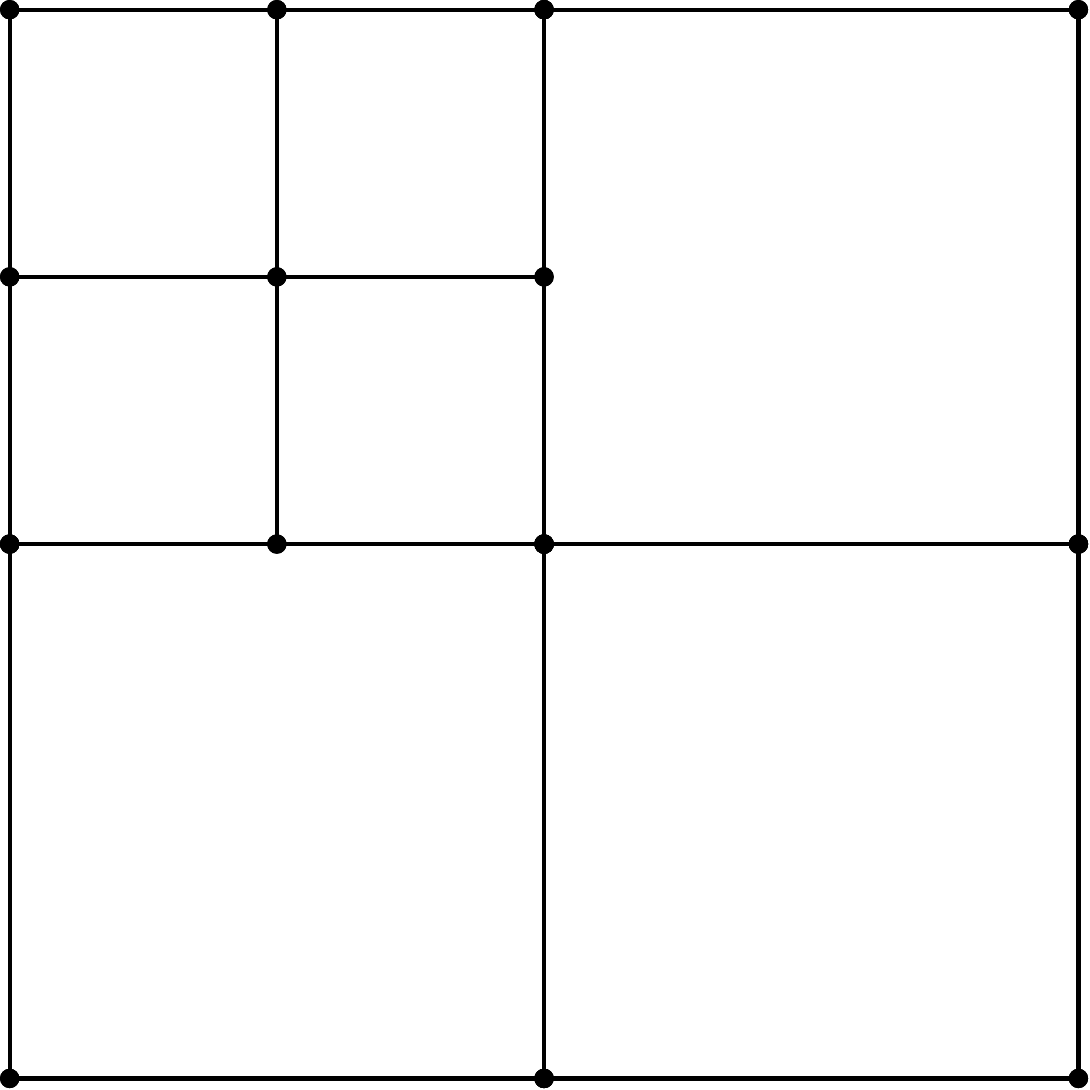}}\hfill
    \subfloat[Octree mesh, SBFEM \label{fig:sketchOct}]{\includegraphics[height=0.23\textwidth]{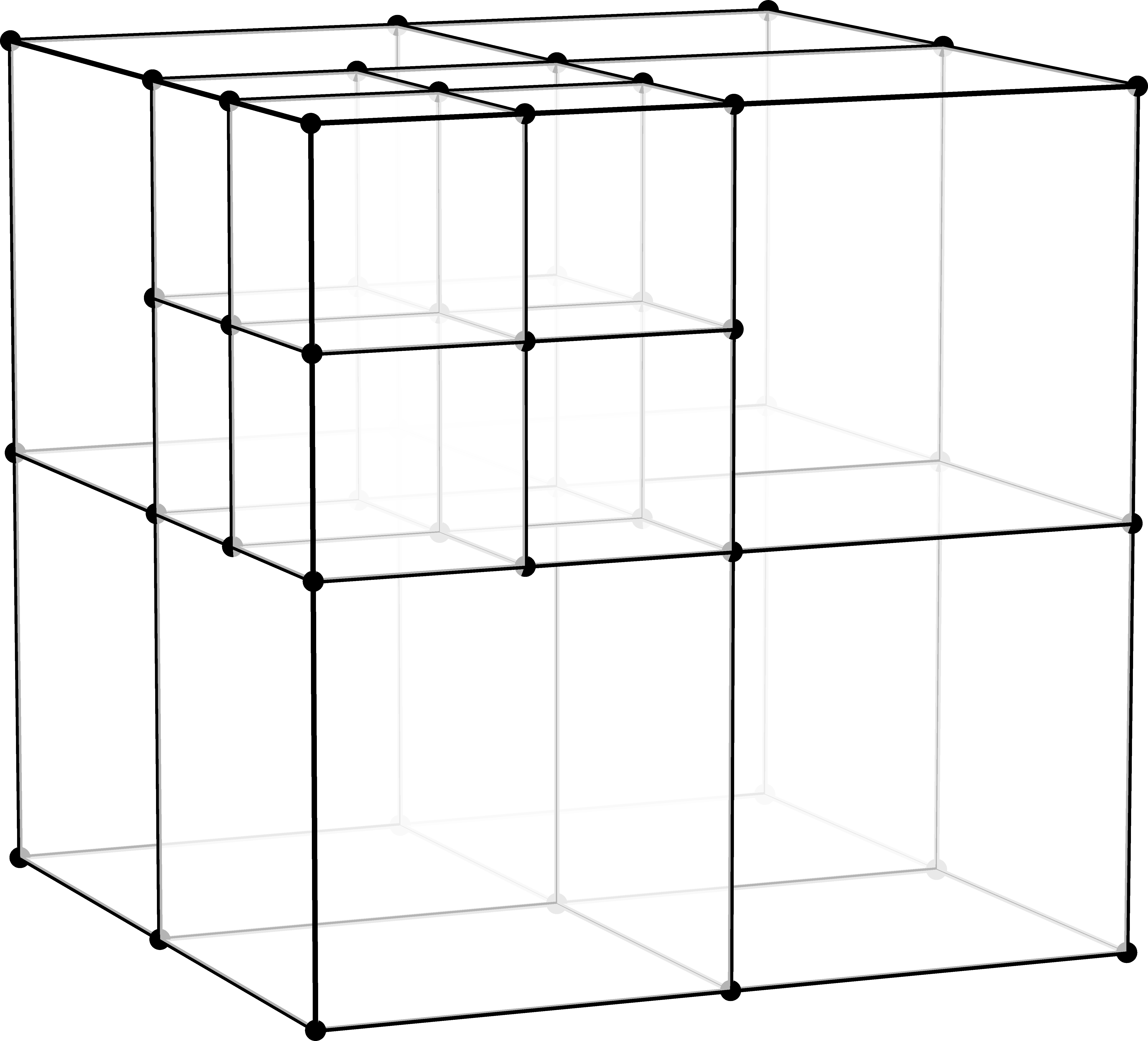}} \hfill
    \subfloat[Quadtree mesh, FEM	\label{fig:sketchQuad_tri}]{\includegraphics[height=0.23\textwidth]{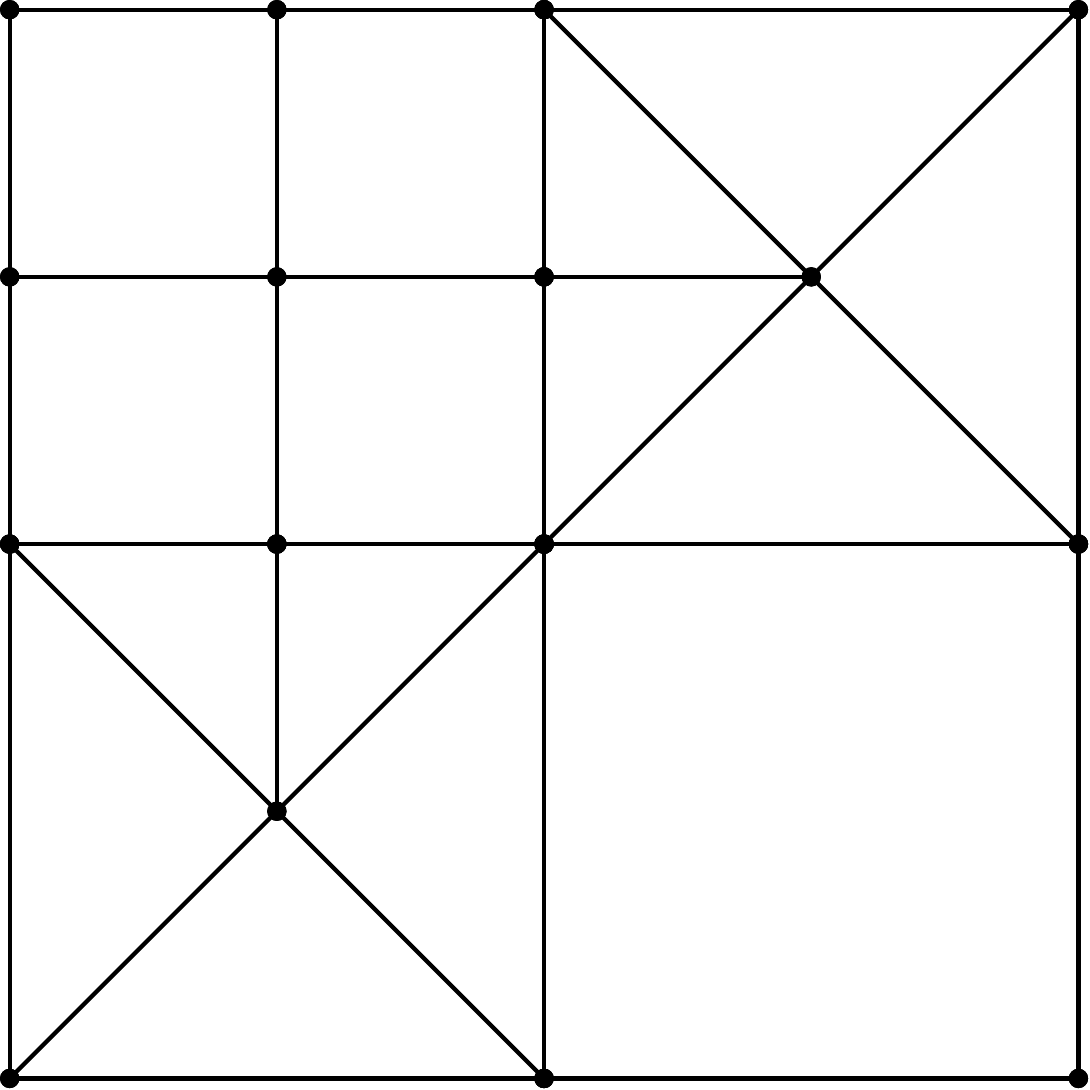}}
    \caption{Discretization in the SBFEM: (a) In 2D, a quadtree-mesh is used to discretize the computational domain straightforwardly. (b) In 3D, surfaces need to be discretized, which feature a structure similar to the quadtree meshes in 2D. (c) In previous work, these surfaces have been triangulated using standard finite elements. Using the transition elements discussed in this paper, the surface meshes can be handled directly without further subdivision (see part (b) of the figure). \label{fig:sketchQuadOct}}
\end{figure}

The SBFEM was originally developed to model large and unbounded domains in the context of soil-structure interaction and was inspired by concepts such as similarity \cite{Wolf1994}, cloning \cite{Wolf1994a}, as well as the thin layer method \cite{Kausel1994b,Kausel1981b}. A detailed description of the underlying formulation can be found in the early papers \cite{Song1997,Wolf2000a,Song2000a}, as well as the recent textbook \cite{Song2018}. Later, it has been noticed that the SBFEM can be employed as a means to construct arbitrary star-convex\footnote{A domain $\Upomega$ is called star-convex if there exists a point $\vt{r}_0\in \Upomega$ such that the line segment from $\vt{r}_0$ to any point in $\Upomega$ is contained in $\Upomega$. In other words, for the domain to be star-convex, there must be a point from where the whole boundary is `visible'. Any convex domain is also star-convex.} elements \cite{Ooi2012a,Ooi2013,Chiong2014,Ooi2014}.  Such polygonal elements enable the application of rather flexible meshing procedures, compared to conventional finite elements which are usually restricted to triangular/quadrilateral and tetrahedral/hexahedral shapes. A particular variant of such a meshing paradigm consists in the use of domain decompositions of the quad-/octree type. The most common of this class of decompositions consists of square-/cube-shaped cells (see Fig.~\ref{fig:sketchQuad} for a minimal example). While, from a geometrical viewpoint, these meshes only require quadrilateral/hexahedral subdomains, the SBFEM concept allows each side to be divided into an arbitrary number of surface elements on the boundary. Hence, coupling subdomains of different size is straightforward. This idea has been used successfully for the efficient meshing of complex geometries, where flexible local refinement is desirable \cite{Man2014,Ooi2015b,Gravenkamp2017a,Gravenkamp2017c}. Such a meshing paradigm is particularly useful for \textit{image-based analyses}, i.e., in applications where the geometry (the distribution of material parameters, etc.) is provided in a pixel graphics format. A number of approaches have been developed to mesh images, most of which can be divided into two categories (see \cite{young2008efficient} for a comprehensive literature review): The first type includes meshing based on boundary detection for each region in the image, where a region is assumed to be homogeneous. This can be straightforwardly achieved using techniques such as a marching cubes algorithm \cite{lorensen1987marching}. Subsequently, these boundaries are utilized to mesh the region inside using conventional meshing algorithms, such as the advancing front \cite{lohner1988generation} or Delaunay tessellation \cite{du2006recent}. The second type involves techniques that mesh the image directly. These techniques are generally quicker due to the straightforward nature of the meshing process. One very basic example is the pixel-based approach \cite{keyak1990automated}, where each pixel is modeled as a quadrilateral/hexahedral finite element. While such a meshing approach is absolutely trivial and fast, it leads to an unnecessarily large number of degrees of freedom. In order to improve the efficiency of the meshing, other direct meshing techniques can be deployed such as the aforementioned quad-/octree meshing structure \cite{young2008efficient}.

Quadtree meshing is performed by recursively dividing an image matrix into four equal-sized cells at a time. A criterion of homogeneity is established based on the difference in the maximum and minimum color intensity within a cell. If this difference exceeds a user-defined threshold, the cell is divided. The maximum and minimum cell size allowed are also specified by the user. As a result, the meshing scheme adaptively refines the regions around different interfaces. Simultaneously, the scheme keeps relatively larger cells within each region where the material can be assumed to be homogeneous.
However, since cells of different sizes exist in quadtree meshes, compatibility issues are encountered when using conventional finite elements. This issue was circumvented by utilizing the scaled boundary finite element method (SBFEM) for quadtree meshes, as this method allows the use of arbitrary polyhedral subdomains. A detailed explanation of image-based analysis using the SBFEM can be found in \cite{Gravenkamp2017a,Saputra2017}.

Notwithstanding the success of this method in two-dimensional image-based analyses, the situation in three-dimensional modeling is somewhat different. The quad\-tree-based domain decomposition can be extended to 3D without much further effort \cite{Saputra2017}. In this case, we divide a cubic cell into 8 smaller ones of equal size and refer to this technique as octree decomposition (Fig.~\ref{fig:sketchOct}). From the viewpoint of the SBFEM, each cubic cell constitutes one subdomain, which implies that we have to discretize the surface of each cube in a finite element sense. Fig.~\ref{fig:sketchOct} reveals that the surfaces of 3D octree decompositions basically consist of 2D quadtree meshes and hence exhibit the same compatibility issues discussed above. That is, due to the two-dimensional finite element grid that has to be used for the discretization of the boundary, hanging nodes are generated. In this example, the mesh on the front surface of the structure in Fig.~\ref{fig:sketchOct} is identical to the two-dimensional case in Fig.~\ref{fig:sketchQuad}. In finite element applications, the hanging node problem can be circumvented by either using some kind of constraints (constraint equations, Lagrange multipliers, etc.) or triangular elements.

In previous works, this issue has been tackled by the second option, i.e., triangulating the surfaces of individual cubes to ensure compatibility with adjacent subdomains, see Fig.~\ref{fig:sketchQuad_tri}. This approach has been applied successfully to solve problems related to acoustics \cite{Liu2019a} and fracture \cite{Saputra2015}. However, this concept may be considered somewhat unsatisfactory as the triangular elements require additional degrees of freedom and are intrinsically less accurate compared to quadrilateral elements. Furthermore, the previous approach does not allow combining elements of different order within the same model, which is one of the major advantages exploited in 2D SBFEM models \cite{Gravenkamp2017a,Gravenkamp2018a}.

To improve on the previous formulation, we present in the current paper an approach to discretize the surfaces of octree-based SBFEM models without introducing additional nodes. This is achieved by employing a certain type of transition element that allows us to formulate shape functions on quadrilateral domains with an arbitrary subdivision of each edge. The application of such elements was inspired by the papers on the so-called `pNh elements' \cite{Graeff-Weinberg1996a,Gabbert1999,Weinberg2002}, however, similar formulations exist and date back to the 1970s \cite{Gordon1973a,ArtcileGordon1973b,ArticleGordon1971,Birkhoff1974,ArticleCavendish1975}. In a recent publication, we presented a detailed review and generalization of these classes of transition elements \cite{Duczek2019b}. There, we refer to these general transition elements as xNy-elements to indicate that it is in principle possible to couple arbitrary element families. In the current paper we shall only provide a very brief summary of the formulation (Sect.~\ref{sec:Transition}) and refer the reader to the pertinent literature. The purpose of this work is mainly to demonstrate that these transition elements can in fact be used to solve the hanging node problem in three-dimensional SBFEM models. The motivation for this approach is twofold:

\begin{itemize}
    \item We wish to avoid any triangulation of the surfaces in three-dimensional (octree-based) models. This goal can easily be achieved using said transition elements, since they allow splitting each edge into an arbitrary number of sections. Consequently, we avoid introducing additional elements that are not required by the topological decomposition. Especially when using higher order interpolation, the reduction in degrees of freedom (DOFs) can be very significant.
    \item We wish to allow different element orders within the same model, which is particularly useful in the case of inhomogeneous materials and wave propagation problems. Up to now, conventional Lagrange elements have been used on the surfaces of SBFEM subdomains, which required the element order to match between adjacent subdomains. Using the proposed approach generally allows connecting elements of arbitrary orders.
\end{itemize}
%
%
\section{Theory}
\label{sec:Theory}
The overarching goal of this contribution is to combine the concepts of SBFEM, octree-based mesh generation, and transition elements. This combination results in an efficient and robust numerical tool which can be applied to a wide range of problems in engineering and physics. To this end, we provide a brief descriptions of the fundamental principles underlying the SBFEM in Sect.~\ref{sec:SBFEM}, while octree-meshes are discussed in Sect.~\ref{sec:Octree}. Considering the formulation of transition elements, we sketch the basic ideas in Sect.~\ref{sec:Transition} and refer to the literature concerning transfinite elements for an in-depth derivation.

\subsection{SBFEM}
\label{sec:SBFEM}
The formulation of the Scaled Boundary Finite Element Method is presented for the case of linear elasticity, governed by the following partial differential equation for the displacement field $\vt{u}=\vt{u}(x,y)$:
\begin{equation}
    \mathcal{L}^\T\m{D}\mathcal{L}\vt{u}(x,y) + \rho \ddot{\vt{u}}(x,y) = \vt{0}
\end{equation}
Here, $\rho$ denotes the mass density and
$\m{D}$ is the elasticity matrix, which, for isotropic material with shear modulus $G$ and Poisson's ratio $\nu$, reads
\begin{equation}
    \mathbf{D}=\frac{2 G}{(1-2 \nu)}\left[\begin{array}{cccccc}
            {1-\nu} & {\nu}   & {\nu}                                                               \\
            {\nu}   & {1-\nu} & {\nu}                                                               \\
            {\nu}   & {\nu}   & {1-\nu}                                                             \\
                    &         &         & \frac{1-2 \nu}{2}                                         \\
                    &         &         &                   & \frac{1-2 \nu}{2}                     \\
                    &         &         &                   &                   & \frac{1-2 \nu}{2} \\
        \end{array}\right]
\end{equation}
$\mathcal{L}$ is the differential operator
\begin{equation}\label{eq:L}
    \mathcal{L}=\left[\begin{array}{cccccc}{\partial_{x}} & {0} & {0} & {\partial_{y}} & {\partial_{z}} & {0} \\ {0} & {\partial_{y}} & {0} & {\partial_{x}} & {0} & {\partial_{z}} \\ {0} & {0} & {\partial_{z}} & {0} & {\partial_{x}} & {\partial_{y}}\end{array}\right]^{\mathrm{T}}
\end{equation}
which may be rewritten for later use as
\begin{equation}
    \mathcal{L}=\hat{\m{b}}_1 \partial_x + \hat{\m{b}}_2 \partial_x + \hat{\m{b}}_3 \partial_z
\end{equation}
with
\begin{align}\nonumber
    \hat{\m{b}}_1 & =\left[\begin{array}{cccccc}
            1 & 0 & 0 & 0 & 0 & 0 \\
            0 & 0 & 0 & 1 & 0 & 0 \\
            0 & 0 & 0 & 0 & 1 & 0
        \end{array}\right]^\T \quad
    \hat{\m{b}}_2=\left[\begin{array}{cccccc}
            0 & 0 & 0 & 1 & 0 & 0 \\
            0 & 1 & 0 & 0 & 0 & 0 \\
            0 & 0 & 0 & 0 & 0 & 1
        \end{array}\right]^\T          \\
    \hat{\m{b}}_1 & =\left[\begin{array}{cccccc}
            0 & 0 & 0 & 0 & 1 & 0 \\
            0 & 0 & 0 & 0 & 0 & 1 \\
            0 & 0 & 1 & 0 & 0 & 0
        \end{array}\right]^\T
\end{align}
By adjusting the differential operator, as well as the material parameters, other common linear partial differential equations are obtained and can be solved by the SBFEM in a similar fashion \cite{Birk2009,Liu2019a}.

The SBFEM constitutes a semi-analytical approach, in the sense that the governing partial differential equation (in three coordinates) is discretized in two directions while remaining analytical in the third coordinate. In order to employ this idea effectively, a particular coordinate transformation is usually applied such that the two discretized directions (denoted as $\eta,\zeta$) define a parametrization of a domain's boundary, while the analytical coordinate $\xi$ describes the direction from the origin to the boundary.\footnote{There are a few exceptions where such a coordinate transformation is not required and the concept of semi-discretization can be applied directly in a Cartesian, polar, or cylindrical coordinate system to model structures of (piece-wise) constant cross-section \cite{Man2012,Gravenkamp2014f,Gravenkamp2014a,Chen2014b,Krome2017,Krome2017a,Gravenkamp2018}} A few key steps of the procedure are summarized in this section for easier reference. The complete formulation including a detailed derivation of the basic coordinate transformation and the semi-discretization can be found in \cite{Song1997}. For dynamic problems, we apply the solution procedure proposed in \cite{Song2009}, which is based on a continued fraction expansion of the dynamic stiffness matrix. A concise summary of the procedures applied to dynamic problems is presented in \cite{Gravenkamp2018a}.

In general, the coordinate transformation employed in the SBFEM is written as
\begin{subequations}
    \begin{align}
        \hat{x}(\xi, \eta, \zeta)=\xi x(\eta, \zeta) \\
        \hat{y}(\xi, \eta, \zeta)=\xi y(\eta, \zeta) \\
        \hat{z}(\xi, \eta, \zeta)=\xi z(\eta, \zeta)
    \end{align}
\end{subequations}
where $\hat{x},\hat{y},\hat{z}$ define a three-dimensional Cartesian coordinate system and $x,y,z$ are the values of these Cartesian coordinates on the boundary of the domain, parametrized by the local coordinates $\eta,\zeta$.\footnote{This slightly peculiar choice of denoting the standard Cartesian coordinates by $\hat{\bullet}$ simply stems from the fact that the coordinates on the boundary occur far more frequently in the formulation of the SBFEM. This notation is adopted here for consistency with \cite{Song1997,Saputra2017}. } Without loss of generality, the origin of the coordinate system is chosen to be inside the domain.  Furthermore, the domain is assumed to be star-convex. The boundary is usually interpolated using a set of two-dimensional shape functions $\mathbf{N}(\eta, \zeta)$ (such as Lagrange shape functions or splines); thus, the coordinates and their derivatives can be written as
\begin{subequations}
    \begin{alignat}{3}
         & \hat{x}(\xi, \eta, \zeta)          &  & =\xi &  & \mathbf{N}(\eta, \zeta)\mathbf{x}_\mathrm{n}          \\
         & \hat{x}_{,\xi}(\xi, \eta, \zeta)   &  & =    &  & \mathbf{N}(\eta, \zeta)\mathbf{x}_\mathrm{n}          \\
         & \hat{x}_{,\eta}(\xi, \eta, \zeta)  &  & =\xi &  & \mathbf{N}_{,\eta}(\eta, \zeta)\mathbf{x}_\mathrm{n}  \\
         & \hat{x}_{,\zeta}(\xi, \eta, \zeta) &  & =\xi &  & \mathbf{N}_{,\zeta}(\eta, \zeta)\mathbf{x}_\mathrm{n}
    \end{alignat}
\end{subequations}
and analogously for $\hat{y},\hat{z}$. The vectors $\mathbf{x}_\mathrm{n},\mathbf{y}_\mathrm{n},\mathbf{z}_\mathrm{n}$ denote the nodal values of the coordinates (or, more generally, the coefficients of shape functions if the chosen basis is not node-based). However, the reader may note that, for the current contribution, an interpolation of the boundary is not required as all surfaces of the octree-meshes addressed here are squares, for which the coordinate transformation could easily be written in closed form. Nevertheless, we stick to this more general formulation for later reference. With these transformations, the Jacobian determinant is obtained as
\begin{equation}
    \hat{J}=\xi J
\end{equation}
where $J$ is the determinant of the matrix
\begin{equation}
    \mathbf{J}=\left[\begin{array}{ccc}
            x         & y          & z          \\
            x_{,\eta} & y_{,\eta}  & z_{,\eta}  \\
            x_{,\xi}  & y_{,\zeta} & z_{,\zeta}
        \end{array}\right]
\end{equation}
The differential operator $ \mathcal{L}$ transformed into the local coordinates reads
\begin{equation}
    \mathcal{L}=\mathbf{b}_1 \partial_\xi+\frac{1}{\xi}
    \left( \mathbf{b}_2 \partial_{\eta} + \mathbf{b}_3 \partial_\zeta \right)
\end{equation}
The displacement field is interpolated on the boundary as
\begin{equation}
    \mathbf{u}(\xi, \eta, \zeta)=\mathbf{N}(\eta, \zeta) \mathbf{u}(\xi)
\end{equation}
where we assume for simplicity that the same shape functions are employed as for the interpolation of the geometry (isoparametric elements). Consequently, the approximated strain-displacement relationship reads
\begin{equation}
    \vg{\varepsilon}=\mathbf{B}_1 \mathbf{u}_\mathrm{n}(\xi)_{,\xi}+\frac{1}{\xi} \mathbf{B}_2 \mathbf{u}_\mathrm{n}(\xi)
\end{equation}
with
\begin{subequations}
    \begin{align}
        \mathbf{B}_1 & =\mathbf{b}_1 \mathbf{N}(\eta, \zeta)                                                       \\
        \mathbf{B}_2 & =\mathbf{b}_2 \mathbf{N}(\eta, \zeta)_{,\eta}+\mathbf{b}_3 \mathbf{N}(\eta, \zeta)_{,\zeta}
    \end{align}
\end{subequations}
Applying the method of weighted residuals to obtain a weak form of the governing equations in the coordinates $\eta,\zeta$, leads to the ordinary matrix differential equation
\begin{multline}\label{eq:ode}
    \mathbf{E}_0 \xi^2 \mathbf{u}_\mathrm{n}(\xi)_{,\xi\xi} + \left((d-1) \mathbf{E}_0-\mathbf{E}_1+\mathbf{E}_1^\mathrm{T}\right) \xi \mathbf{u}_\mathrm{n}(\xi)_{,\xi}\\ +
    \left((d-2) \mathbf{E}_1^{\mathrm{T}}-\mathbf{E}_2\right) \mathbf{u}_\mathrm{n}(\xi) + \omega^2\xi^2\mathbf{M}_0 \mathbf{u}_\mathrm{n}(\xi)=0
\end{multline}
with the coefficient matrices
\begin{subequations}
    \begin{align}
        \mathbf{E}_0   & =\int_S \mathbf{B}_1^\mathrm{T} \mathbf{D} \mathbf{B}_1 J \mathrm{d} \eta \mathrm{d} \zeta         \\
        \mathbf{E}_1   & =\int_S \mathbf{B}_{2}^{\mathrm{T}} \mathbf{D} \mathbf{B}_{1} J \mathrm{d} \eta \mathrm{d} \zeta   \\
        \mathbf{E}_{2} & =\int_{S} \mathbf{B}_{2}^{\mathrm{T}} \mathbf{D} \mathbf{B}_{2} J \mathrm{d} \eta \mathrm{d} \zeta \\
        \mathbf{M}_0   & =\int_S \mathbf{N}^\mathrm{T} \rho \mathbf{N} J \mathrm{d} \eta \mathrm{d} \zeta
    \end{align}
\end{subequations}
where $d$ denotes the number of dimensions (2 or 3). Many details on the properties and the available solution procedures of similar differential equations are discussed in \cite{Kausel2019}. For the static case ($\omega = 0$), the differential equation is usually re-written as
\begin{equation}\label{eq:odeZ}
    \xi \mathbf{X}(\xi)_{, \xi}=-\mathbf{Z X}(\xi)
\end{equation}
with
\begin{equation}
    \mathbf{Z}=\left[\begin{array}{cc}{\mathbf{E}_{0}^{-1} \mathbf{E}_{1}^{\mathrm{T}}-0.5(d-2) \mathbf{I}} & {-\mathbf{E}_{0}^{-1}} \\ {-\mathbf{E}_{2}+\mathbf{E}_{1} \mathbf{E}_{0}^{-1} \mathbf{E}_{1}^{\mathrm{T}}} & {-\left(\mathbf{E}_{1} \mathbf{E}_{0}^{-1}-0.5(d-2) \mathbf{I}\right)}\end{array}\right]
\end{equation}
and
\begin{equation}
    \mathbf{X}(\xi)=
    \left[\begin{array}{c}
            \xi^{+0.5(d-2)} \mathbf{u}_\mathrm{n}(\xi) \\
            \xi^{-0.5(d-2)} \mathbf{q}_\mathrm{n}(\xi)
        \end{array}\right]
\end{equation}
Here, $\mathbf{q}_\mathrm{n}$ denotes the vector of nodal forces and $\mathbf{I}$ is the identity matrix. A set of eigenfunctions of the semi-discrete differential operator in Eq.~\eqref{eq:odeZ} can be obtained by performing an eigenvalue decomposition of the matrix $\mathbf{Z}$. Alternatively -- to improve accuracy in cases where (near)-parallel eigenvectors exist -- a block-diagonal Schur factorization can be employed instead, which leads to a decomposition of the form
\begin{equation}
    \mathbf{Z}\left[\begin{array}{ll}{\mathbf{\Psi}_{11}} & {\mathbf{\Psi}_{12}} \\ \mathbf{\Psi}_{21} & {\mathbf{\Psi}_{22}}\end{array}\right] =
    \left[\begin{array}{ll}{\mathbf{\Psi}_{11}} & {\mathbf{\Psi}_{12}} \\ {\mathbf{\Psi}_{21}} & {\mathbf{\Psi}_{22}}\end{array}\right]\left[\begin{array}{cc}{\mathbf{S}_{11}} & {0} \\ {0} & {\mathbf{S}_{22}}\end{array}\right]
\end{equation}
Details are discussed in \cite{Song2004}. In the above decomposition, $\mathbf{\Psi}_{11}, \mathbf{\Psi}_{21}$ denote the modes that lead to finite displacements at the origin, and $\boldsymbol{S}_{11}$ is the corresponding Schur block. The solution for the case of a bounded domain is then obtained as
\begin{subequations}
    \begin{align}
        \mathbf{u}(\xi) & =\mathbf{\Psi}_{11} \xi^{-\left(\mathbf{S}_{11}+0.5(d-2) \mathbf{I}\right)} \mathbf{c}_{1} \\ \mathbf{q}(\xi) &=\mathbf{\Psi}_{21} \xi^{-\left(\mathbf{S}_{11}-0.5(d-2) \mathbf{I}\right)} \mathbf{c}_{1}
    \end{align}
\end{subequations}
On the boundary ($\xi=1$), we obtain
\begin{subequations}
    \begin{align}
        \mathbf{u}_\mathrm{b} & =\mathbf{\Psi}_{11} \mathbf{c}_1 \\
        \mathbf{q}_\mathrm{b} & =\mathbf{\Psi}_{21} \mathbf{c}_1
    \end{align}
\end{subequations}
which we use to compute a stiffness matrix
\begin{equation}
    \mathbf{K}=\mathbf{\Psi}_{21} \mathbf{\Psi}_{11}^{-1}
\end{equation}
such that
\begin{equation}
    \mathbf{q}_\mathrm{b}=\mathbf{K} \mathbf{u}_\mathrm{b}
\end{equation}
In the dynamic case ($\omega > 0$), the differential equation in displacements \eqref{eq:ode} is converted into one for the dynamic stiffness $\mathbf{S}(\omega)$ on the boundary:
\begin{multline}\label{eq:dynstiff}
    (\m{S}(\omega)- \Eo)\Ez^{-1}(\m{S}(\omega)-\Eo^\T)- \Et \\
    + (d-2)\m{S}(\omega) \omega \m{S}(\omega)_{,\omega}
    +\omega^2 \M=0
\end{multline}
Assuming an approximation of the form
\begin{equation}\label{eq:Sapprox}
    \m{S}(\omega)\approx \K -\omega^2 \m{M}
\end{equation}
a mass matrix is obtained from the solution of the Lyapunov equation:
\begin{equation}
    (\K-\Eo)\Ez^{-1}\m{M}+\m{M}\Ez^{-1}(\K-\Eo^\T)+d\cdot\m{M}-\M=0
\end{equation}
To enhance the accuracy for higher frequencies, high-order stiffness and mass matrices can be computed of the form
\begin{subequations}\label{eq:eqmot_time_def}
    \begin{align}
        \m{K}_\mathrm{h} & =\mathrm{diag}(\m{K},\m{S}_0^{(1)},\m{S}_0^{(2)},...,\m{S}_0^{(M_\mathrm{cf})}]) \\
        \m{M}_\mathrm{h} & =\left[
            \begin{array}{ccccc}
                \m{M}             & -\m{X}^{(1)}      & 0             & \cdots & 0                         \\
                -[\m{X}^{(1)}]^\T & \m{S}_1^{(1)}     & -\m{X}^{(2)}  & \cdots & 0                         \\
                0                 & -[\m{X}^{(2)}]^\T & \m{S}_1^{(2)} & \cdots & 0                         \\
                \vdots            & \vdots            & \vdots        & \ddots & \vdots                    \\
                0                 & 0                 & 0             & \cdots & \m{S}_1^{(M_\mathrm{cf})}
            \end{array}\right]
    \end{align}
\end{subequations}
where the matrices $\m{S}_0^{(i)},\m{S}_1^{(i)}$ are derived from a continued fraction expansion of the dynamic stiffness matrix
\begin{multline}\label{eq:S_cf_impr}
    \m{S}(\omega)=\K -\omega^2 \m{M}-\omega^4 \m{X}^{(1)} (\m{S}_0^{(1)}-\omega^2\m{S}_1^{(1)}\\
    -\omega^4\m{X}^{(2)}(\m{S}_0^{(2)}-\omega^2 \m{S}_1^{(2)}-...
    -\omega^4 \m{X}^{(M_\mathrm{cf})}(\m{S}_0^{(M_\mathrm{cf})}\\
    -\omega^2 \m{S}_1^{(M_\mathrm{cf})})^{-1}[\m{X}^{(M_\mathrm{cf})}]^\T )^{-1}[\m{X}^{(2)}]^\T )^{-1}[\m{X}^{(1)}]^\T
\end{multline}
and the matrices $\m{X}^{(i)}$ are introduced for preconditioning. The details of this approach can be found in \cite{Song2009,Birk2012}.

\subsection{Octree mesh}
\label{sec:Octree}
While the concept of the SBFEM can generally be applied to any star-convex subdomain, the use of the quadtree/octree decomposition has become particularly popular. As already discussed in the introduction, octree-meshes can be obtained very easily by starting with a coarse discretization consisting of few cubes and recursively dividing cubes into eight smaller cubes according to some criteria. These criteria may be based on whether the cube is intersected by a boundary or whether the variation of material parameters within the cube exceeds a predefined threshold. For further details, we refer again to \cite{Saputra2017}. Also, the reader may note that such a meshing procedure will not lead to smooth boundaries of complex geometries due to the `stair-case'-like approximation of interfaces. These geometry errors are acceptable in many applications -- particularly in image-based analyses where the available initial geometry description is a pixel-based image. In other cases, it is desirable to improve the description of the boundaries by applying smoothing algorithms, cutting the cubes that are intersected by the boundary \cite{Liu2017}, or combining the octree decomposition with fictitious domain concepts \cite{Gravenkamp2017c}. Such approaches are the subject of other publications and will not be discussed here in detail. In the current contribution, we restrict ourselves to pure octree-meshes since this is the part where our approach differs from previous works.

\begin{figure*}[h]
    \centering
    \subfloat{\includegraphics[width=0.15\textwidth]{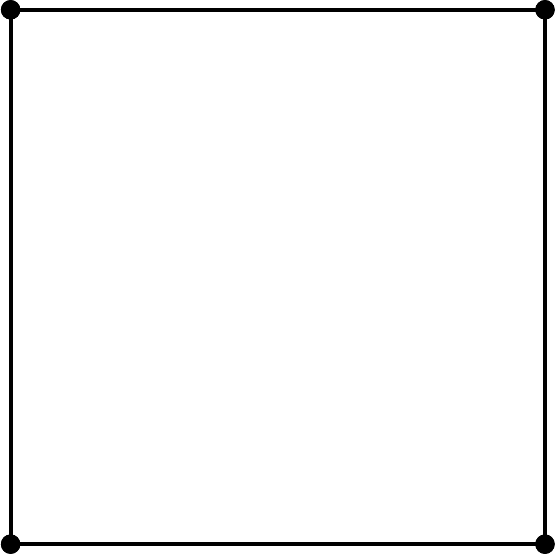}}\hfill
    \subfloat{\includegraphics[width=0.15\textwidth]{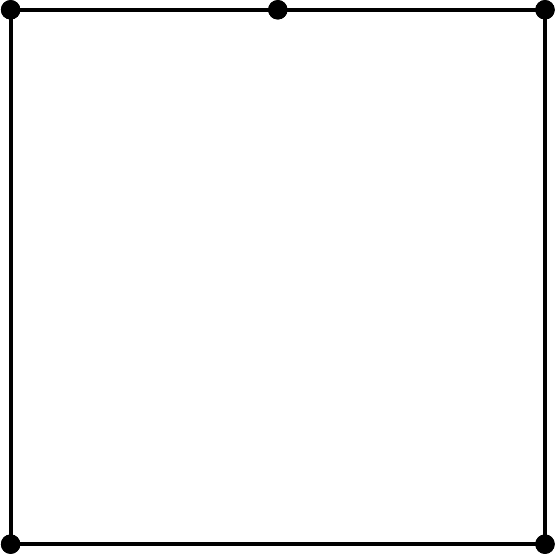}}\hfill
    \subfloat{\includegraphics[width=0.15\textwidth]{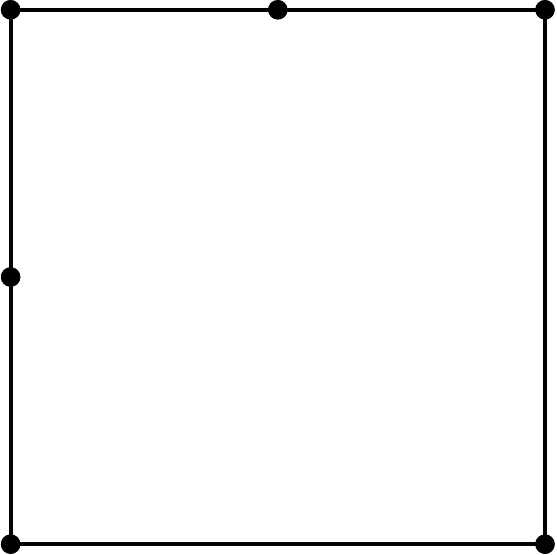}}\hfill
    \subfloat{\includegraphics[width=0.15\textwidth]{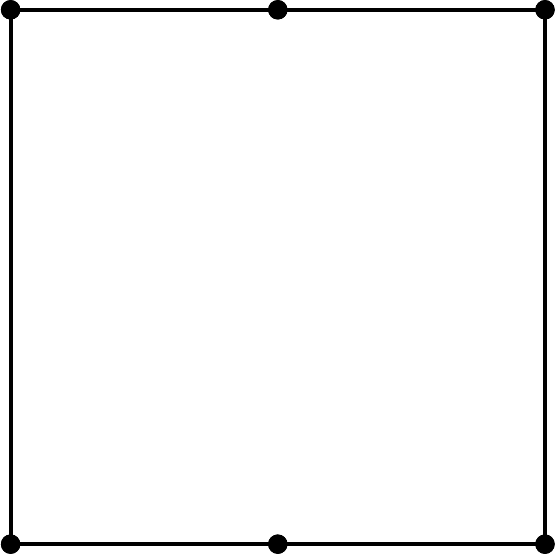}}\hfill
    \subfloat{\includegraphics[width=0.15\textwidth]{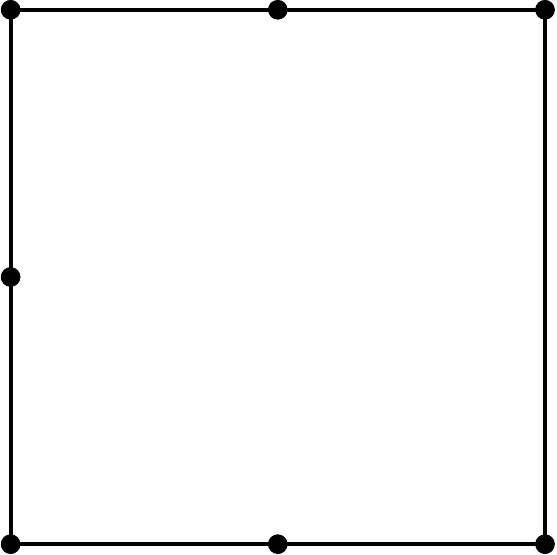}}\hfill
    \subfloat{\includegraphics[width=0.15\textwidth]{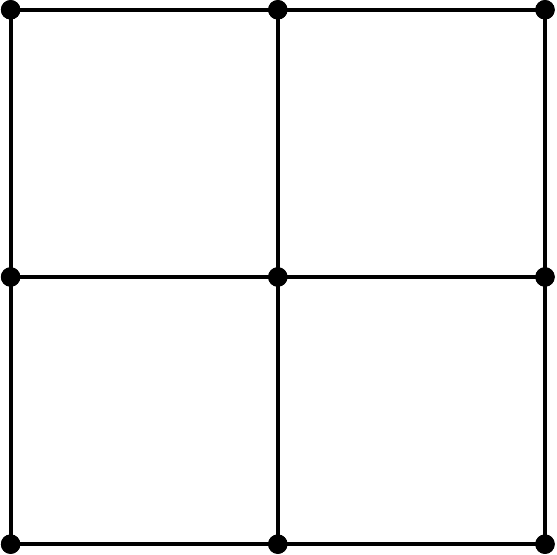}}\\
    \subfloat{\includegraphics[width=0.15\textwidth]{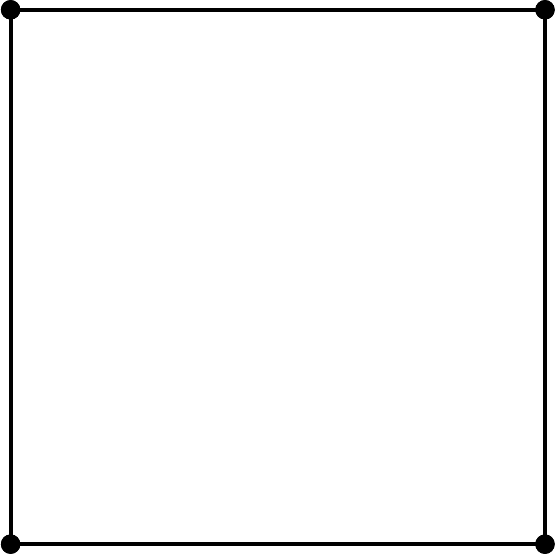}}\hfill
    \subfloat{\includegraphics[width=0.15\textwidth]{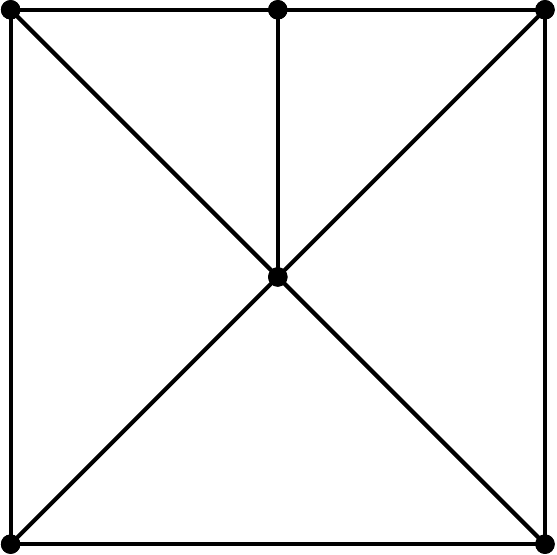}}\hfill
    \subfloat{\includegraphics[width=0.15\textwidth]{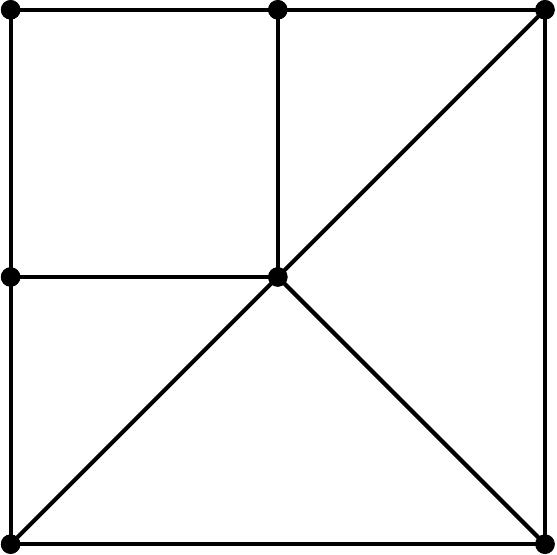}}\hfill
    \subfloat{\includegraphics[width=0.15\textwidth]{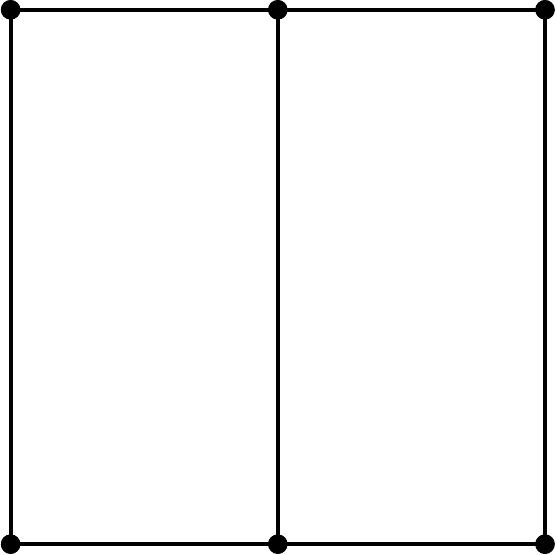}}\hfill
    \subfloat{\includegraphics[width=0.15\textwidth]{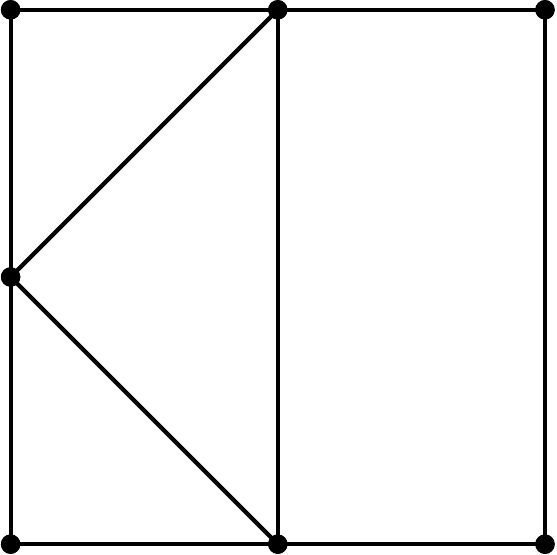}}\hfill
    \subfloat{\includegraphics[width=0.15\textwidth]{quadtreePatternTri6_p1}}
    \caption{Mesh patterns on the surfaces of a balanced octree decomposition. Transition elements allow us to use quadrilateral elements in all cases (top), while conventional elements require a subdivision of the surfaces to treat hanging nodes (bottom). \label{fig:QuadtreePatternsP1}}
\end{figure*}

The proposed application of transition elements changes the way the meshes on the surfaces of the SBFEM subdomains are handled. For now, let us assume that the octree decomposition is `balanced', which means that adjacent subdomains differ by at most one refinement level. Consequently, each edge of each cube is either split into two segments or not. Taking into account the symmetries under rotation, there exist only six distinct meshing patterns on the cube surfaces. These patterns are depicted in Fig.~\ref{fig:QuadtreePatternsP1} and differ in the number of sides that are split in order to ensure conformal coupling to the adjacent subdomain. Note that the first and last pattern in each row all consist of standard quadrilateral elements; thus, in the proposed approach, there are only five different types of surface elements to be considered.\footnote{To be more precise, there are two different cases in which each of the four edges of a given element is split into two segments: If smaller elements need to be connected only to the edges, these edges need to be segmented, while an additional node on the surface is not actually required. If, however, the surface itself is connected to smaller cubes, we need to split this surface element into four. In our current implementation, these two cases are treated the same way by splitting the surface into four elements.} As will be shown in the following sections, we can construct shape functions with nodes at the desired positions. For comparison, Fig.~\ref{fig:QuadtreePatternsP1} also displays the subdivision of the patterns using conventional triangular and quadrilateral elements as has been done in previous works.

\begin{figure*}[htb]
    \centering
    \subfloat{\includegraphics[width=0.15\textwidth]{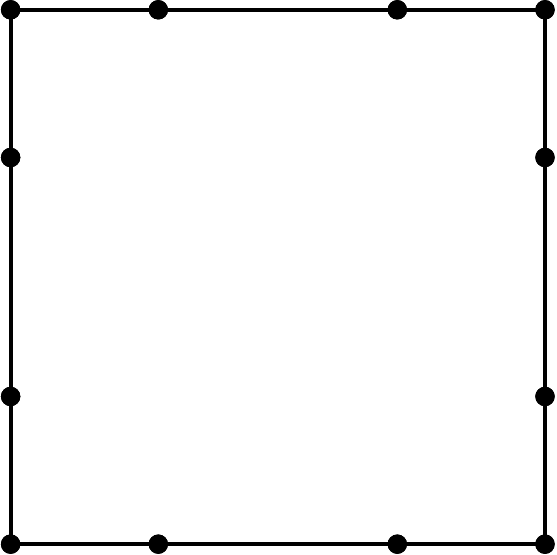}}\hfill
    \subfloat{\includegraphics[width=0.15\textwidth]{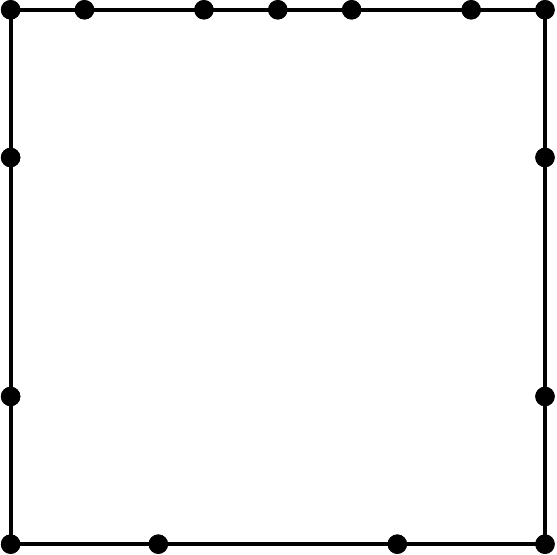}}\hfill
    \subfloat{\includegraphics[width=0.15\textwidth]{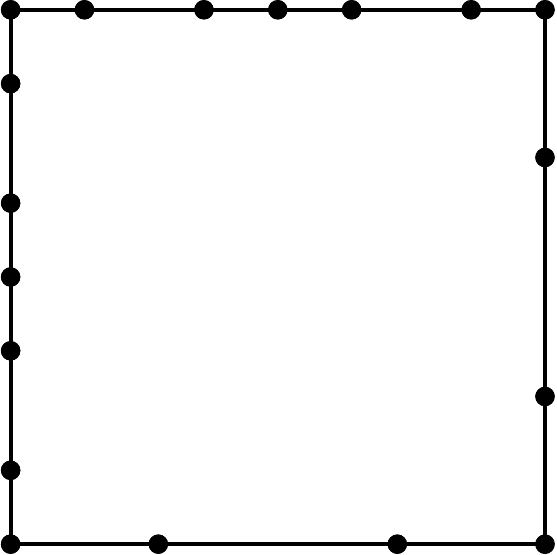}}\hfill
    \subfloat{\includegraphics[width=0.15\textwidth]{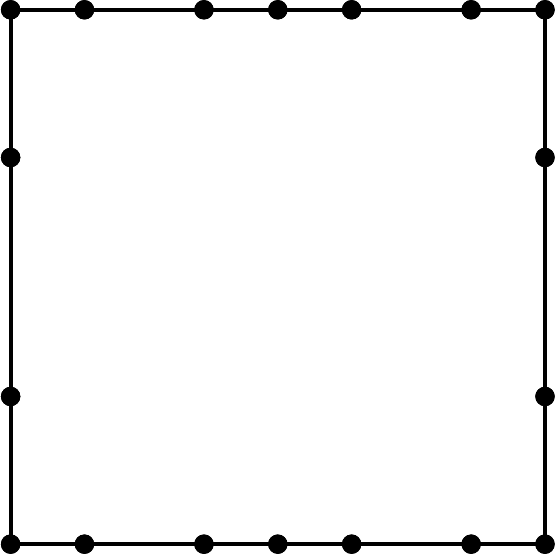}}\hfill
    \subfloat{\includegraphics[width=0.15\textwidth]{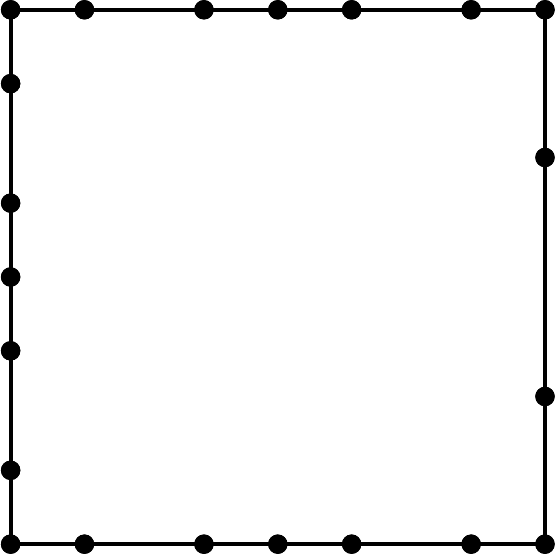}}\hfill
    \subfloat{\includegraphics[width=0.15\textwidth]{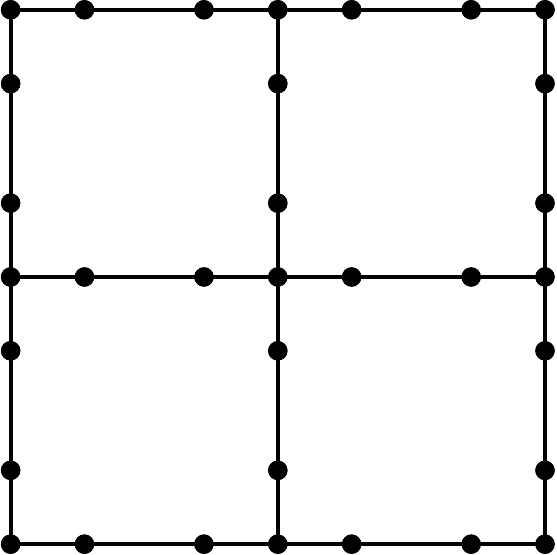}}\\
    \subfloat{\includegraphics[width=0.15\textwidth]{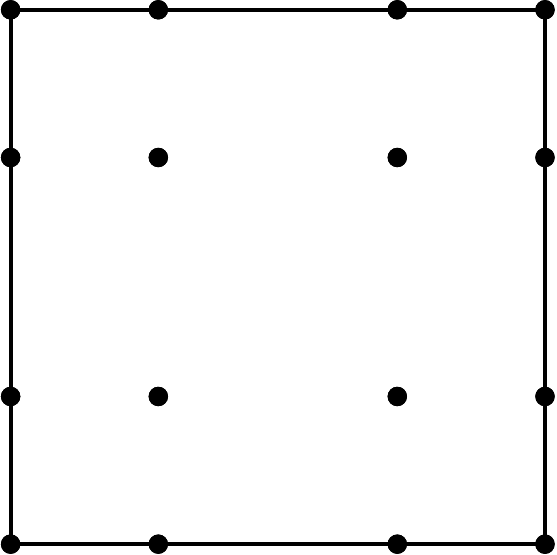}}\hfill
    \subfloat{\includegraphics[width=0.15\textwidth]{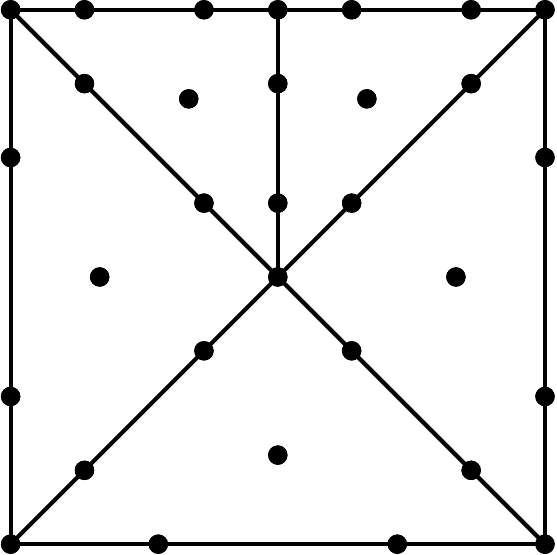}}\hfill
    \subfloat{\includegraphics[width=0.15\textwidth]{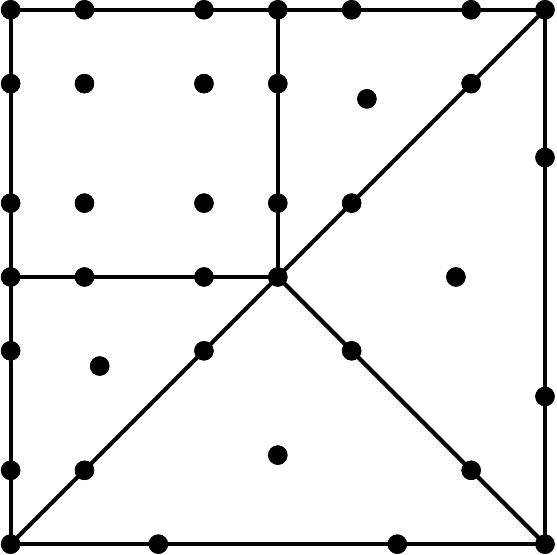}}\hfill
    \subfloat{\includegraphics[width=0.15\textwidth]{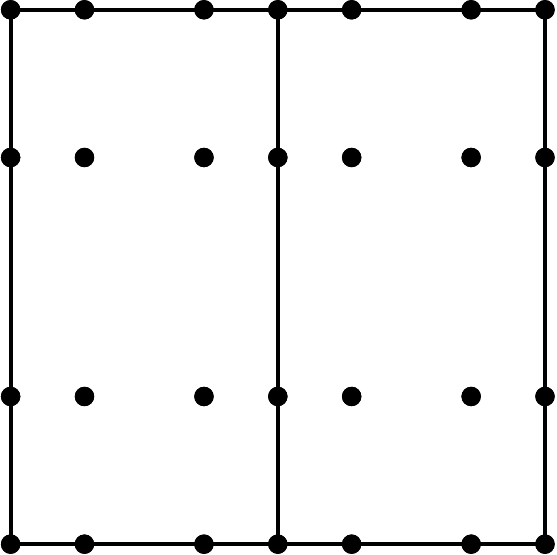}}\hfill
    \subfloat{\includegraphics[width=0.15\textwidth]{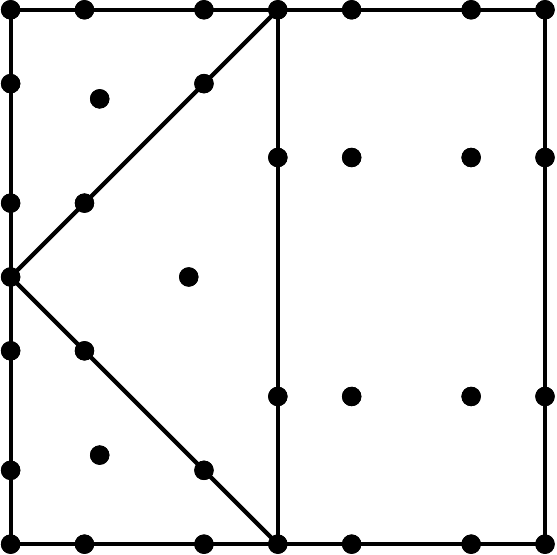}}\hfill
    \subfloat{\includegraphics[width=0.15\textwidth]{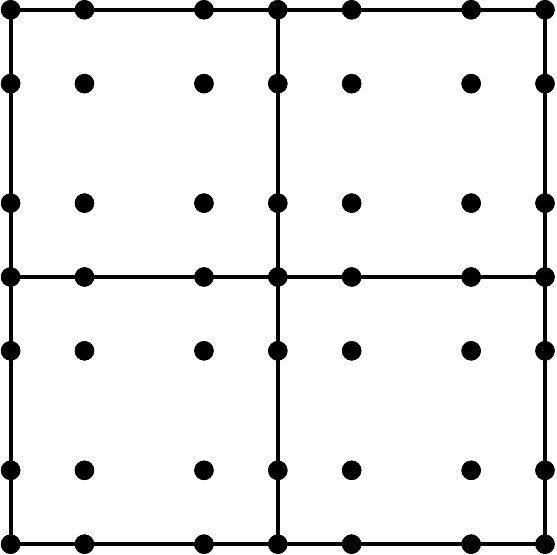}}
    \caption{Mesh patterns for the case of third order interpolation: Transition elements (top) and Lagrange elements (bottom). \label{fig:QuadtreePatternsP3}}
\end{figure*}

The difference in both approaches becomes prominent when applying higher-order interpolations. As an example, Fig.~\ref{fig:QuadtreePatternsP3} shows the same patterns with additional nodes using a polynomial order of 3. Note that the transition elements do not require internal nodes to achieve complete polynomials up to an order of 3. This is in contrast to the Lagrange elements employed previously. In that sense, the shape functions are similar to that of serendipity finite elements. Furthermore, it should also be noted that, when using transition elements, there is no need to restrict the decomposition to balanced meshes. In fact, subdomains of arbitrary sizes can be coupled straightforwardly. This would be rather tedious when using conventional elements, since adequate triangulations would have to be defined for each of the (arbitrarily large number of) different patterns.

\subsection{Transition elements based on the \emph{xNy} element concept}
\label{sec:Transition}
As already discussed in the previous sections, transition elements are being used in the context of finite element models to realize local mesh refinement. In the framework of the SBFEM, transition elements can be used to discretize each face of a hexahedral SBFE which is connected to a number of smaller SBFEs having the same or similar shape functions. This versatility would, in principle, allow us to connect different types of SBFEs based on various sets of shape functions, comprehensively discussed in Ref.~\cite{Gravenkamp2019}. In this paper, our focus is, however, on high order shape functions based on Lagrange polynomials that are defined on a non-equidistant grid of points \cite{BookKarniadakis2005,BookPozrikidis2014}. Typically, a Gauss-Lobatto-Legendre (GLL) nodal distribution is chosen to generate the high order Lagrange shape functions. This particular set of nodes includes the element (interval) boundaries at $\pm 1$, while the interior nodes are the roots of the Lobatto polynomials of order $p\,{-}\,1$ if shape functions of order $p$ are constructed \cite{PhDDuczek2014}.

\subsubsection{Projection: Fundamental idea}
\label{sec:ProjectionIdea}
The idea of transition elements, as implemented in this paper, can be traced back to the works of Gordon and co-workers~\cite{ArticleGordon1971,ArtcileGordon1973a,ArtcileGordon1973b,BookGordon1982}. In their articles, the transition element was referred to as transfinite element due to the special type of mapping that was applied. In the remainder of the article, we will stick to the term transition element. It is worth mentioning that the same approach to deriving shape functions for transition elements can also be used to achieve an accurate approximation of the geometry in high order FEMs. In the \emph{p}-FEM, this kind of geometry description is commonly referred to as blending function method \cite{BookSzabo1991,ArticleKiralyfalvi1997,PhDBroeker2001} and is widely used to incorporate the exact geometry of the structure stemming, for example, from computer-aided-design (CAD) software. More recent applications of the transition element technology can be found in Refs.~\cite{PhDWeinberg1996,ArticleWeinberg2002,ArticleProvatidis2006,ArticleProvatidis2011,BookProvatidis2019}

In the current contribution, the main goal is to couple SBFEs of different sizes, where each boundary of each element is discretized using two-dimensional quadrilateral elements. Therefore, the task is to develop a transition element that features piece-wise polynomial shape functions such that it can couple to multiple smaller elements. At this point, only the basic idea is sketched and we refer to the pertinent literature~\cite{ArticleGordon1971,ArtcileGordon1973a,ArtcileGordon1973b,BookGordon1982,PhDWeinberg1996,ArticleWeinberg2002,ArticleProvatidis2006,ArticleProvatidis2011,BookProvatidis2019} and the references cited therein. In order to achieve conforming elements, we need to interpolate a function $\Xi(\xi,\eta)$ over the reference domain $\Omega$:$\, [-1,1] \times [-1,1]$. To this end, we project the arbitrary bivariate function $\Xi(\xi,\eta)$ onto a different (carefully selected) space of bivariate functions \cite{ArtcileGordon1973b}. Such a projection will be denoted as $\mathcal{P}[\Xi(\xi,\eta)]$ and is composed of two projectors $\mathcal{P}_{\xi}[\Xi(\xi,\eta)]$, $\mathcal{P}_{\eta}[\Xi(\xi,\eta)]$ that interpolate $\Xi(\xi,\eta)$ along the local directions $\xi$ and $\eta$. Additionally, a mixed projector $\mathcal{P}_{\xi\eta}[\Xi(\xi,\eta)]$ is introduced to ensure that redundant terms cancel out. Hence, the projection is written as
\begin{equation}
    \mathcal{P}[\Xi(\xi,\eta)] = \mathcal{P}_{\xi}[\Xi(\xi,\eta)] + \mathcal{P}_{\eta}[\Xi(\xi,\eta)] - \mathcal{P}_{\xi\eta}[\Xi(\xi,\eta)]\,,
    \label{eq:Projector}
\end{equation}
where the mixed projection operator is defined as
\begin{equation}
    \mathcal{P}_{\xi\eta}[\Xi(\xi,\eta)] = \mathcal{P}_{\xi}[\mathcal{P}_{\eta}[\Xi(\xi,\eta)]]\,.
    \label{eq:MixedProjector}
\end{equation}
Note that the individual projection operators are both linear
\begin{equation}
    \mathcal{P}_s[f(\xi,\eta) + g(\xi,\eta)] = \mathcal{P}_s[f(\xi,\eta)] + \mathcal{P}_s[g(\xi,\eta)]
    \label{eq:ProjectorLinear}
\end{equation}
and idempotent
\begin{equation}
    \mathcal{P}_s[\mathcal{P}_s[f(\xi,\eta)]] = \mathcal{P}_s[f(\xi,\eta)]\,,
    \label{eq:ProjectorIdempotent}
\end{equation}
where $f(\xi,\eta)$ and $g(\xi,\eta)$ are (continuous) bivariate functions and the subscript $s\in \{\xi,\eta\}$.

\subsubsection{Projection: Definition of operators}
\label{sec:ProjectionOperator}
The methodology sketched in the previous subsection is deployed to construct conformal finite element spaces ensuring that elements of different sizes can be coupled by introducing an approach to derive piece-wise polynomial shape functions. To achieve this goal, a simple definition of the projection operators suffices. As already introduced in the early works of Gordon and Hall \cite{ArtcileGordon1973a,ArtcileGordon1973b}, a linear Lagrange interpolation polynomial is appropriate to define the projectors as
\begin{align}
    \mathcal{P}_{\xi}[\Xi(\xi,\eta)]  & = \psi_1(\xi) \Xi(\xi_1,\eta) + \psi_2(\xi) \Xi(\xi_2,\eta)\,, \label{eq:ProjectorLagrangeXi}    \\
    \mathcal{P}_{\eta}[\Xi(\xi,\eta)] & = \psi_1(\eta) \Xi(\xi,\eta_1) + \psi_2(\eta) \Xi(\xi,\eta_2)\,, \label{eq:ProjectorLagrangeEta}
\end{align}
where $\xi_i$ and $\eta_i$ denote the nodes that are used to define the Lagrange polynomials. The functions $\psi_i(s)$ are identical to the linear Lagrange polynomials $N_i^1(s)$ and in this context they are referred to as linear blending functions
\begin{align}
    \label{eq:linearBlending1}
    \psi_1(s) & = N_1^1(s) = \cfrac{1}{2}(1-s)\,, \\
    \label{eq:linearBlending2}
    \psi_2(s) & = N_2^1(s) = \cfrac{1}{2}(1+s)\,.
\end{align}
The definition of the mixed projection operator immediately follows from the successive application of Eqs.~\eqref{eq:ProjectorLagrangeXi} and \eqref{eq:ProjectorLagrangeEta} to the function $\Xi(\xi,\eta_i)$
\begin{equation}
    \begin{split}
        & \mathcal{P}_{\xi\eta}[\Xi(\xi,\eta)] =\\ \; & \psi_1(\xi)\psi_1(\eta) \Xi(\xi_1,\eta_1) + \psi_2(\xi)\psi_1(\eta) \Xi(\xi_2,\eta_1)\; + \\
        & \psi_2(\xi)\psi_2(\eta) \Xi(\xi_2,\eta_2) + \psi_1(\xi)\psi_2(\eta) \Xi(\xi_1,\eta_2)\,.
    \end{split}
    \label{eq:ProjectorLagrangeProduct}
\end{equation}
Looking at Eq.~\eqref{eq:ProjectorLagrangeProduct}, it is immediately clear that  the mixed projector interpolates the function only at discrete nodal locations, whereas the individual projectors take the element edges into account. This procedure ensures that elements are compatible at the edges. Based on the presented concept, we are now able to derive the shape functions for arbitrary transition elements. To this end, the function $\Xi(\xi,\eta)$ is chosen such that it represents the piece-wise polynomial function that enables coupling elements of different sizes. Recall that the goal is to couple three-dimensional SBFE-subdomains of different sizes. In this case, it suffices to ensure the edges feature conformal shape functions as long as serendipity elements of degree $p\leq 3$ are used. These elements do not include interior/bubble shape functions that are necessary to the polynomial completeness \cite{ArticleDuczek2019a}. If higher order elements or tensor product formulations \cite{ArticleDuczek2019b} are deployed, we also have to account for the interior/bubble shape functions to ensure a conformal coupling of three-dimensional SBFE-subdomains. Note that the shape functions connected to the interior of the surface elements do not need to be included in the projection process, i.e., the bubble functions are identical to those of standard finite elements. One possibility is to only add those hierarchic bubble modes that are known from the \emph{p}-version of FEM \cite{BookSzabo1991} to complete the polynomial.
%
%
\begin{figure}[ht]
    \centering
    \includegraphics[width=0.4\textwidth,trim=0cm -2cm 0cm 0cm, clip ]{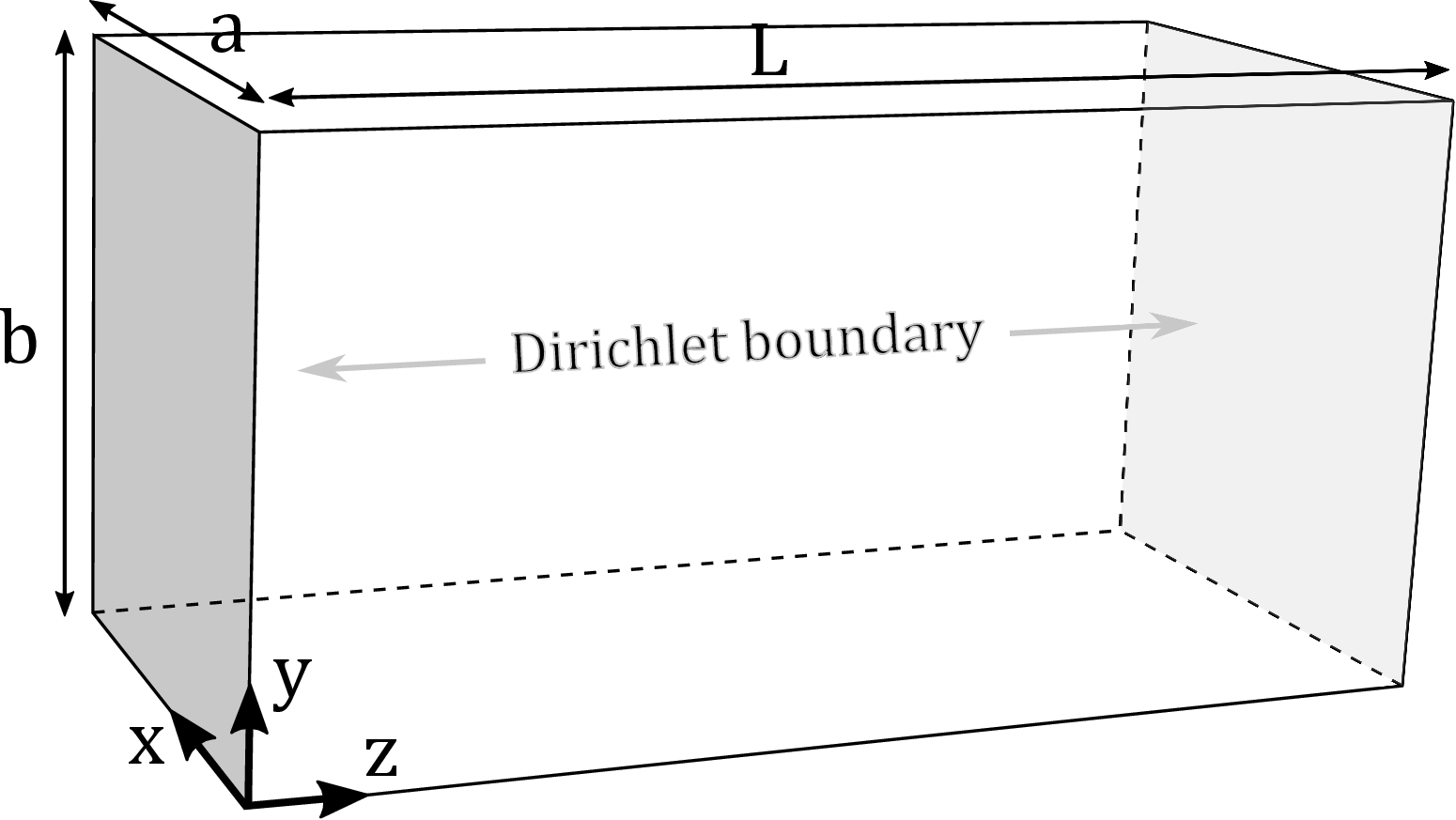}
    \caption{Geometry used to perform patch tests.\label{fig:patchtests_geo}}
\end{figure}

\section{Numerical examples}
\label{sec:NumExamples}
In this section, we present the results of several numerical experiments that we performed in order to validate the proposed approach. We begin by analyzing simple geometries -- namely a rectangular cuboid and a cube -- and conduct standard patch tests, as well as a modal analysis. We then proceed to demonstrate the applicability of the approach to more complex geometries. In all examples, we assume linear elastic behavior of the material and small deformations. We denote by $\Upomega$ the computational domain and $\Upgamma_\mathrm{u}$, $\Upgamma_\mathrm{q}$ are the parts of the boundary where Dirichlet and Neumann boundary conditions are applied, respectively. Hence, the general problem statement may be written as
\begin{alignat}{3}
    \mathcal{L}^\T \mathbf{D} \mathcal{L} \boldsymbol{u} - \rho \ddot{\boldsymbol{u}}- \boldsymbol{f} & =\boldsymbol{0}                                   &  & \qquad \text{in\ } \Upomega            \\
    \boldsymbol{u}                                                                                    & = \boldsymbol{u}_\Upgamma                         &  & \qquad \text{on\ } \Upgamma_\mathrm{u} \\
    \boldsymbol{n}^\T \boldsymbol{\sigma}                                                             & =  \boldsymbol{n}^\T \boldsymbol{\sigma}_\Upgamma &  & \qquad \text{on\ } \Upgamma_\mathrm{q}
\end{alignat}
where $\boldsymbol{u}_\Upgamma$ and $\boldsymbol{\sigma}_\Upgamma$ denote the boundary conditions and $\mathcal{L}$ is the differential operator as given in Eq.~\eqref{eq:L}.

\begin{figure*}[htb]
    \centering
    \subfloat[1$^\text{st}$ order patch test -- uniaxial tension.\label{fig:patchtests_uni}]{\includegraphics[width=0.4\textwidth]{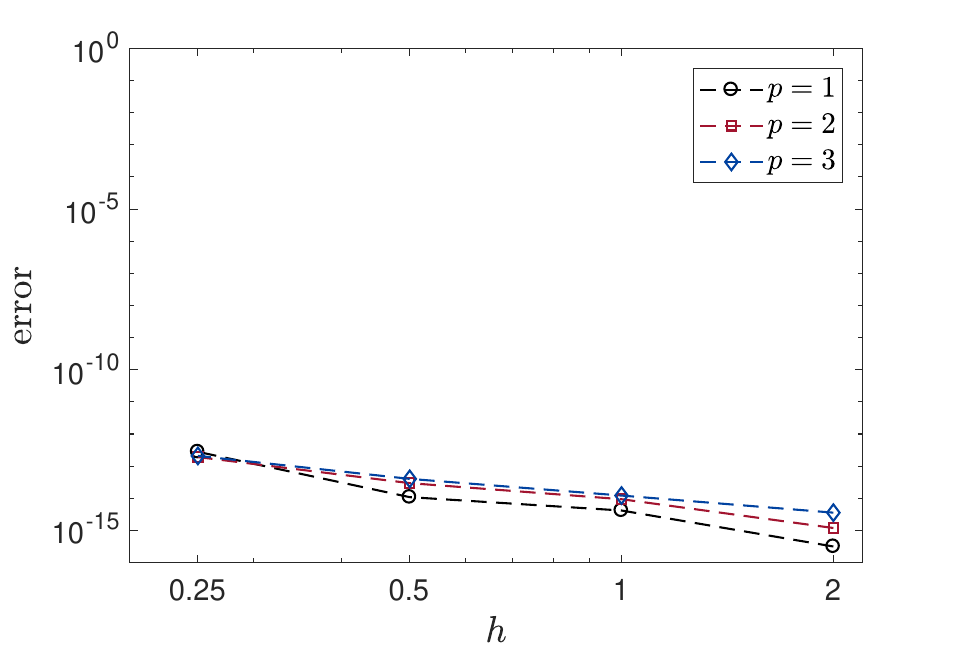}}\qquad
    \subfloat[Exemplary mesh, $h=2$, $p=1$. \label{fig:patchtests_mesh1}]{\includegraphics[width=0.4\textwidth,trim=0cm -2cm 0cm 0cm, clip ]{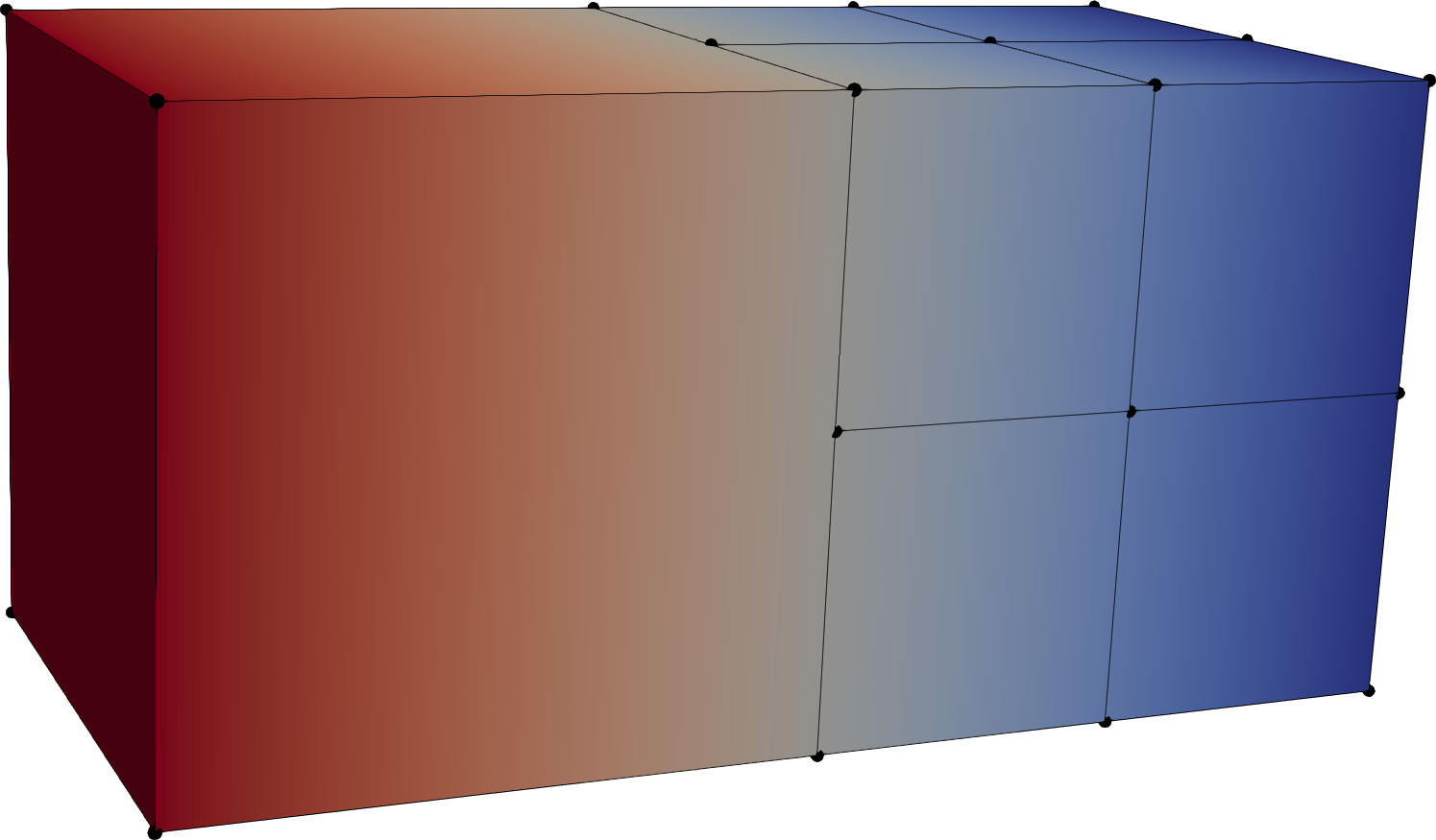}}\\
    \subfloat[2$^\text{nd}$ order patch test -- bending.\label{fig:patchtests_bend}]{\includegraphics[width=0.4\textwidth]{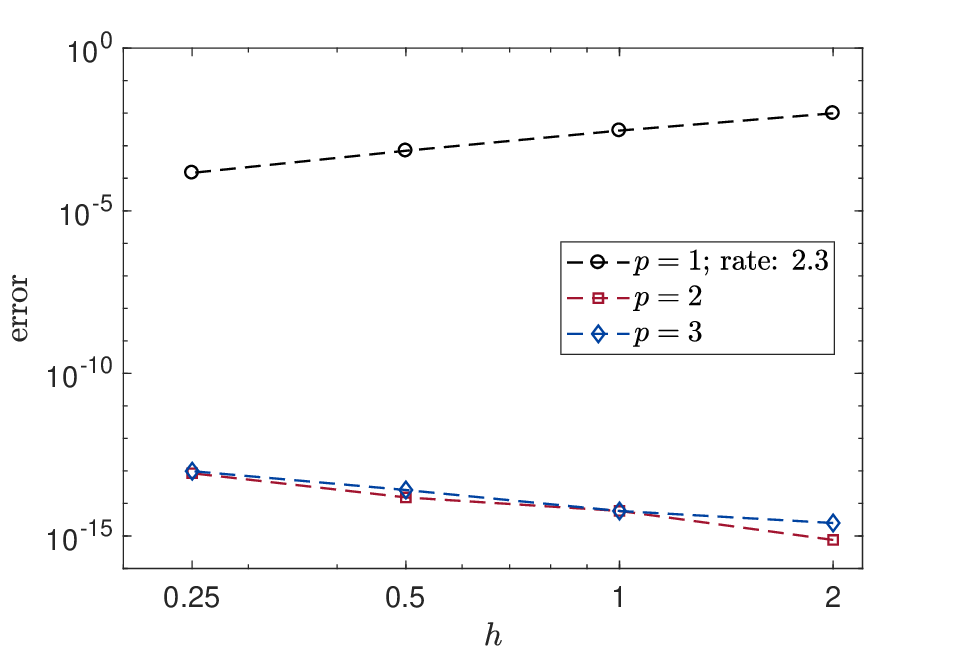}}\qquad
    \subfloat[Exemplary mesh, $h=1$, $p=2$. \label{fig:patchtests_mesh2}]{\includegraphics[width=0.4\textwidth,trim=0cm -2cm 0cm 0cm, clip]{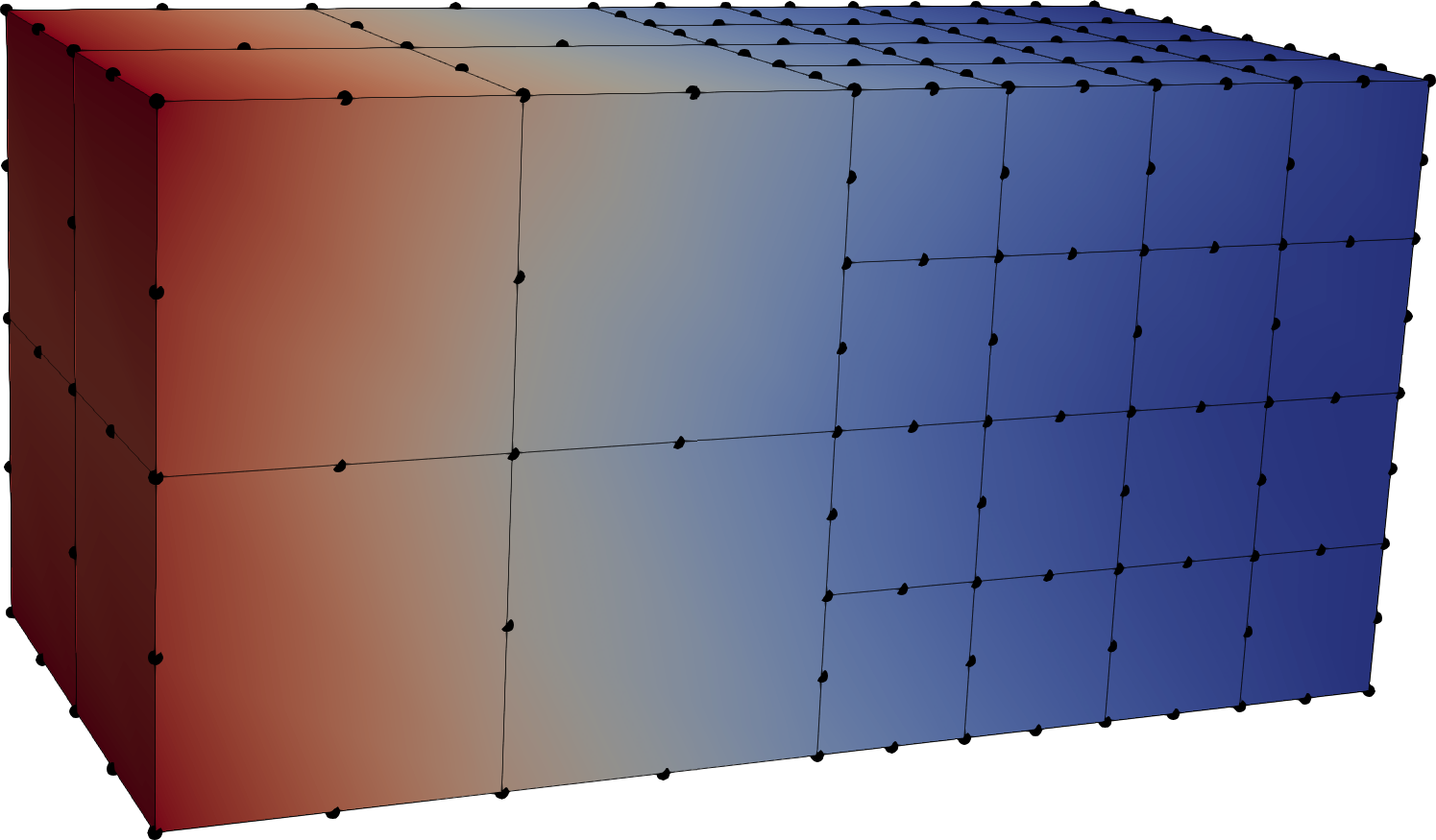}}\\
    \subfloat[3$^\text{rd}$ order patch test -- cantilever beam. \label{fig:patchtests_cant}]{\includegraphics[width=0.4\textwidth]{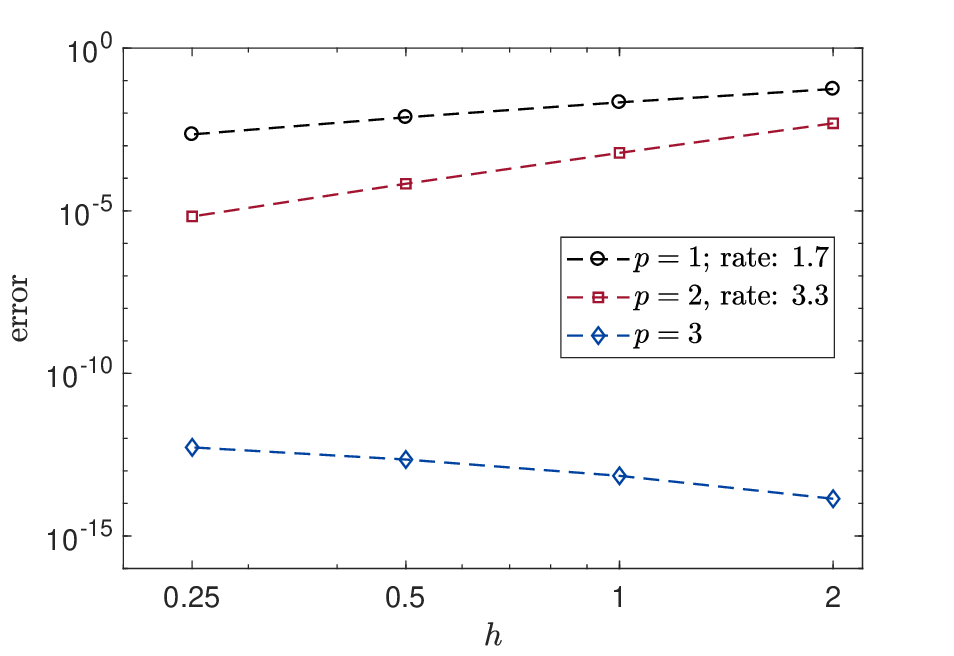}}\qquad
    \subfloat[Exemplary mesh, $h=0.5$, $p=3$. \label{fig:patchtests_mesh3}]{\includegraphics[width=0.4\textwidth,trim=0cm -2cm 0cm 0cm, clip]{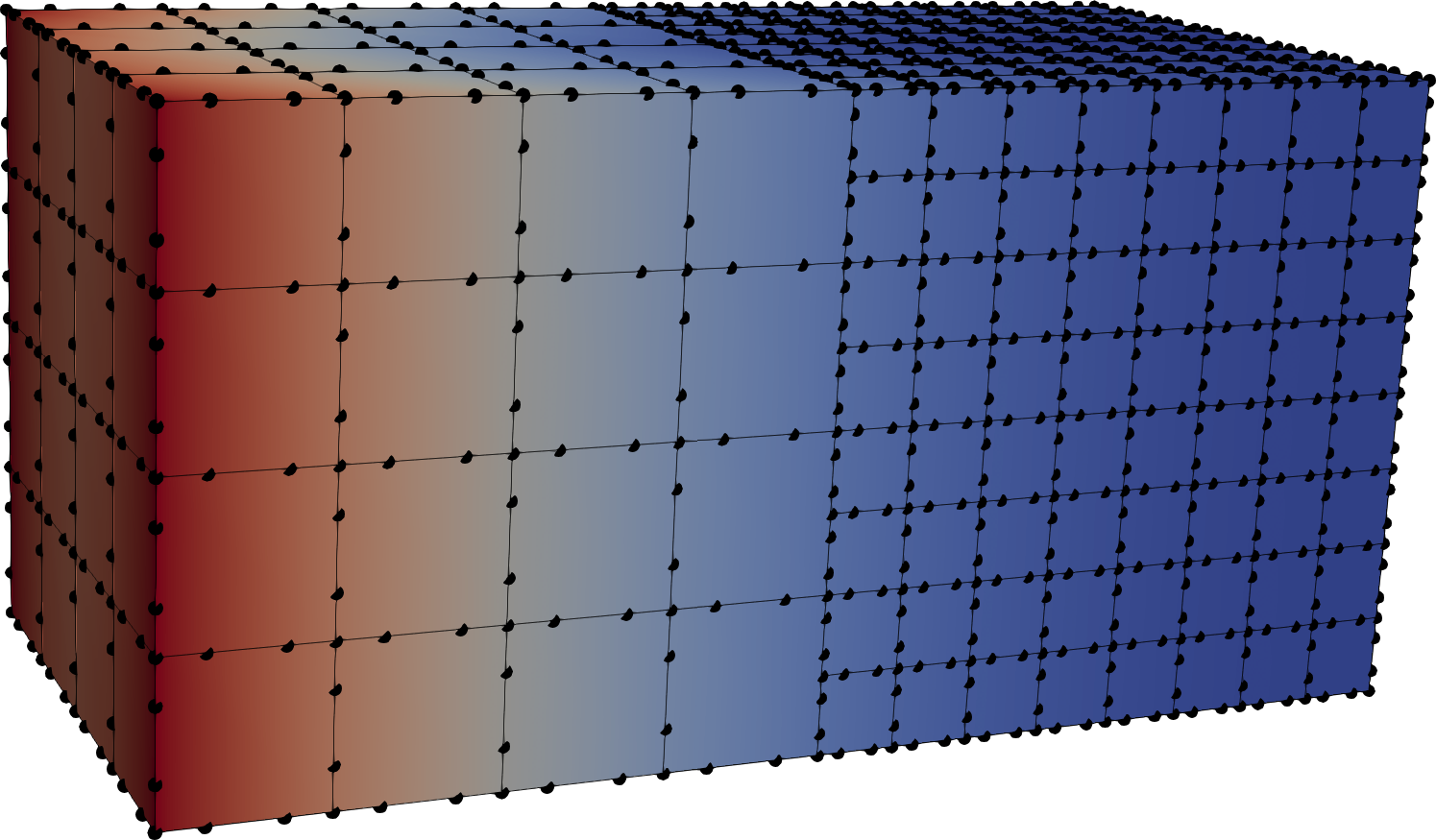}}
    \caption{Results of the 1$^\text{st}$, 2$^\text{nd}$, and 3$^\text{rd}$ order patch tests.\label{fig:patchtests}}
\end{figure*}

\subsection{Static analyses and patch tests}
We begin by conducting linear and higher order patch tests to check the validity of the proposed approach and evaluate the rate of convergence under optimal conditions. These tests are performed for a cuboid with a width and height of $a=b=2$ and a length of $L=4$, see Fig.~\ref{fig:patchtests_geo}. Hence, the computational domain is given as
\begin{equation}
    \Upomega=\left\{\left. (x, y, z) \in \mathbb{R}^3\,\right|\,0\leq x\leq 2,\,0\leq y\leq 2,\,0\leq z\leq 4\right\} \nonumber
\end{equation}
In all tests, Dirichlet boundary conditions are applied to the two faces indicated in Fig.~\ref{fig:patchtests_geo}, while all other surfaces are traction-free. The material is isotropic and homogeneous with the following properties:\\

\begin{tabular}{lc}
    \text{Young's modulus:} & $E =1$   \\
    \text{Poisson's ratio:} & $\nu =0$
\end{tabular}\\

\noindent
Despite the fact that the material is homogeneous, we choose to divide the domain into two regions, where different element sizes are employed. Thus, transition elements are incorporated at the interface between the two regions (see Fig.~\ref{fig:patchtests} for examples). The element order is varied between $p=1$ and $p=3$. We perform $h$-refinement by dividing each subdomain into eight in each refinement step. The element size $h$ is defined as the edge length of the largest element in the mesh. To assess the accuracy of the computed results, we evaluate the $L_2$ norm of the relative error in displacements with respect to an analytical solution. We study the following three cases:

\paragraph{Uniaxial tension}
We assume a state of uniaxial tension, such that the exact solution of this problem is
\begin{equation}
    \boldsymbol{u}^\mathrm{ref}=(0,0,z)^\T
\end{equation}
and consequently, the following Dirichlet boundary conditions are applied:
\begin{align}
    \boldsymbol{u}_{\Upgamma_1} = \boldsymbol{u}^\mathrm{ref}(x,y,0) & = (0,0,0)^\T \\
    \boldsymbol{u}_{\Upgamma_2} = \boldsymbol{u}^\mathrm{ref}(x,y,L) & = (0,0,L)^\T
\end{align}
Hence, the numerical errors should be negligible, as long as the elements on the surface of each subdomain are capable of representing a linear variation of the displacement field exactly. Figure~\ref{fig:patchtests_uni} shows the computed errors when performing $h$-refinement with varying element order. The results show that the proposed approach passes the linear patch test with error levels in the order of $10^{-13}$.

\paragraph{Beam bending} As a 2$^\text{nd}$ order patch test, we analyze the cuboid (with the same properties as before) in a state of pure bending (with a unit radius of curvature). The analytical solution of the displacement field is given in \cite{Timoshenko1951} (Chapter 10) as
\begin{align}
    u_x^\mathrm{ref} & = -\frac{1}{2} ( z^2+\nu (x^2-y^2) ) \\
    u_y^\mathrm{ref} & = -xy \nu                            \\
    u_z^\mathrm{ref} & =\phantom{+}xz
\end{align}
Hence, for the case of $\nu=0$, our Dirichlet boundary conditions read
\begin{align}
    \boldsymbol{u}_{\Upgamma_1} & = \vt{u}^\mathrm{ref}(x,y,0) = \left(0,0,0\right)                 \\
    \boldsymbol{u}_{\Upgamma_2} & = \vt{u}^\mathrm{ref}(x,y,L) = \left(-\frac{1}{2} L^2,0,xL\right)
\end{align}
Results are presented in Fig.~\ref{fig:patchtests_bend}. Since the exact solution is a quadratic function of $x,y,z$, the errors are negligible when using an element order of $p\geq 2$. On the other hand, if an element order of $p=1$ is chosen, the results are not exact. In this case, the numerical convergence rate -- computed based on the last two points of the graphs -- is obtained with a value of 2.3.

\paragraph{Cantilever beam}
Finally, we perform a 3$^\text{rd}$ order patch test by applying boundary conditions that result in a cubic variation of the displacement field. The analytical reference solution of the displacement field $\vt{u}^\mathrm{ref}(x,y,z) = (u_x^\mathrm{ref},u_y^\mathrm{ref},u_z^\mathrm{ref})$ is based on a cantilever beam which is weakly fixed on one end and subject to a shear force $F$ applied at the other end (at $z=4$ in our example). The analytical solution can be found in \cite{Bishop2014}:
\begin{align}
     & u_x^\mathrm{ref}=-\frac{F \nu}{E I} x y z                                                      \\
     & u_y^\mathrm{ref}=\frac{F}{E I}\left[\frac{\nu}{2}\left(x^2-y^2\right) z-\frac{1}{6} z^3\right]
\end{align}
\begin{multline}
    u_z^\mathrm{ref}=\frac{F}{E I} \left[\frac{1}{2} y\left(\nu x^2+z^2\right)+\frac{1}{6} \nu y^3 \right.\\
        \left.+(1+\nu)\left(b^2 y-\frac{1}{3} y^3\right)-\frac{1}{3} a^2 \nu y  \right.\\
        - \left.\frac{4 a^3 \nu}{\pi^3} \sum_{n=1}^\infty \frac{(-1)^n}{n^3} \cos \left(n \pi x / a\right) \frac{\sinh \left(n \pi y / a\right)}{\cosh (n \pi b / a)} \right]
\end{multline}
According to this analytical solution, Dirichlet boundary conditions are applied at the two ends of the cantilever as
\begin{align}
    \boldsymbol{u}_{\Upgamma_1} & =
    \frac{F}{EI}\left(0,\; 0,\; b^2 y-\frac{1}{3} y^3\right)^\T                                                                          \\
    \boldsymbol{u}_{\Upgamma_2} & = \frac{F}{EI}\left(0,\; -\frac{1}{6}L^3,\; -\frac{1}{3}y^3+\left[\frac{1}{2}L^2+b^2\right]y\right)^\T
\end{align}
Again, all other surfaces are traction-free.
As expected, the error is negligible when using an element order of $p=3$, see Fig.~\ref{fig:patchtests_cant}. For linear and quadratic elements, we obtain a numerical convergence rate of $1.7$ and $3.3$, respectively.

\begin{figure}
    \centering
    \subfloat[Exemplary mesh with triangulation. \label{fig:patchtest_triVsPnhb}]{\includegraphics[width=0.4\textwidth,trim=0cm -2cm 0cm 0cm, clip ]{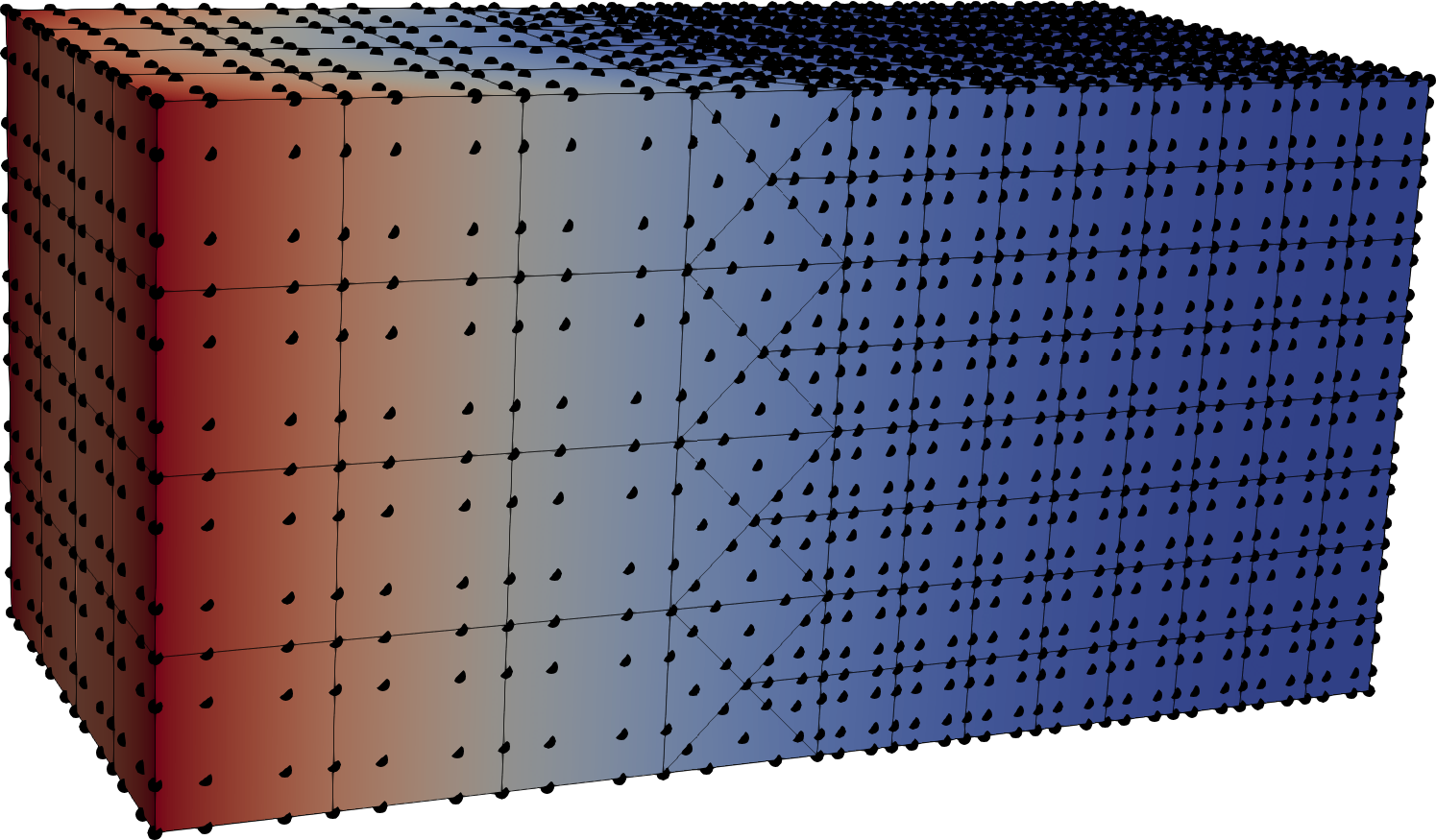}}\\
    \subfloat[Error vs. element size\label{fig:patchtest_triVsPnha}]{\includegraphics[width=0.4\textwidth]{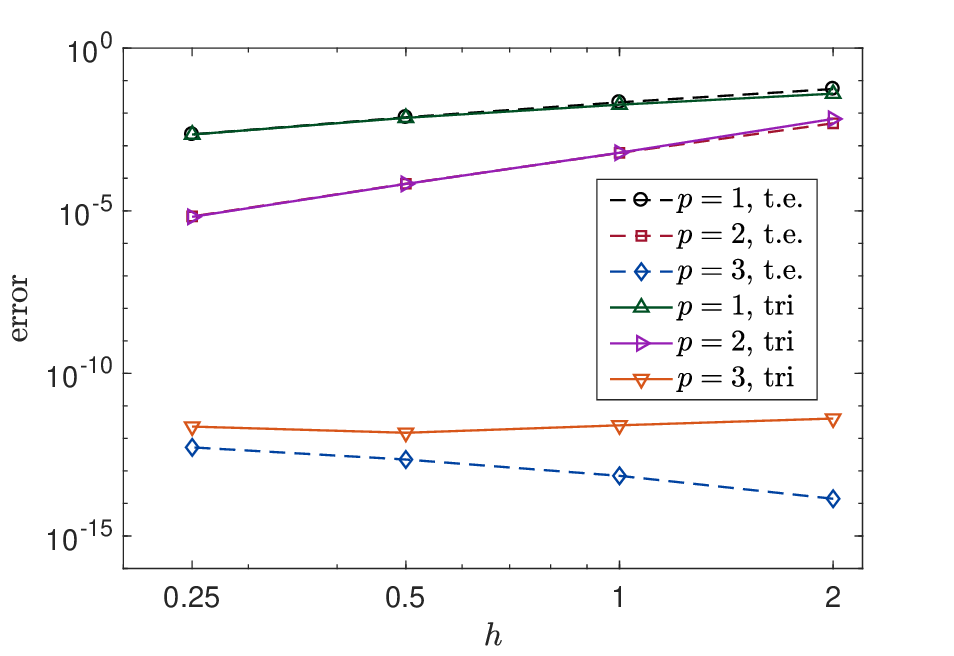}}\\
    \subfloat[Error vs. degrees of freedom.\label{fig:patchtest_triVsPnhc}]{\includegraphics[width=0.4\textwidth]{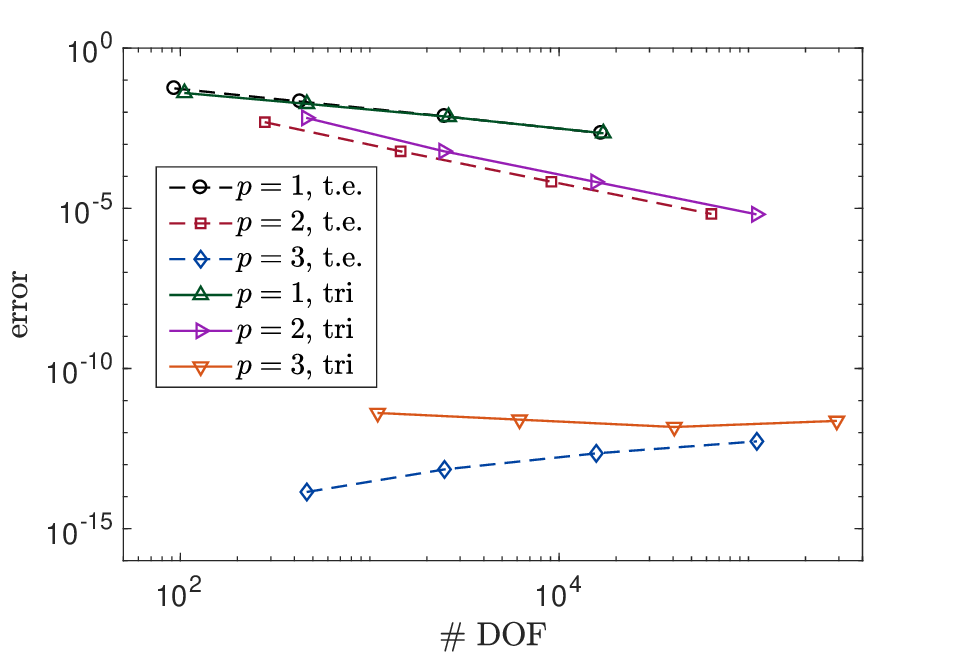}}
    \caption{3$^\text{rd}$ order patch tests -- comparison between transition elements (`xNy') and triangulated surfaces (`tri').\label{fig:patchtest_triVsPnh}}
\end{figure}
\paragraph{Comparison against triangulation}
In addition to the above validation against analytical solutions, we compare the proposed approach to a previous implementation, in which surfaces are subdivided into triangular and rectangular elements rather than using transition elements (see the explanations in Sect.~\ref{sec:Octree}). Such a discretization is depicted in  Fig.~\ref{fig:patchtest_triVsPnhb} and can be compared directly with the corresponding mesh in Fig.~\ref{fig:patchtests_mesh3}. For the sake of brevity, we present here only the results for the cantilever beam. A comparison of the error in displacements using both approaches is presented in Fig.~\ref{fig:patchtest_triVsPnha}. The results obtained based on the proposed discretization (denoted as `xNy' in the figure) are of course identical to the ones in Fig.~\ref{fig:patchtests_cant} and are included here to facilitate the comparison. It can be seen that the application of triangulated surfaces (`tri')\footnote{For conciseness, we refer to these meshes as `triangulated' even though some of the surfaces are subdivided into rectangles to create conforming meshes.} results in very similar error levels for the same element size. Only when using an element order of $p=3$ we observe a significant increase of the error levels -- particularly for large element sizes -- when using triangular elements. Plotting the same errors against the number of degrees of freedom (Fig.~\ref{fig:patchtest_triVsPnh}), the triangulated mesh obviously requires more DOFs to achieve the same error level as the transition elements. This is due to the fact that the quadtree patterns lead to a larger number of interior nodes when applying the triangulation (cf.~Figs.~\ref{fig:QuadtreePatternsP1} and \ref{fig:QuadtreePatternsP3}). While in the case of linear elements (in the current example) only a few additional nodes are required at the centroids of triangulated surfaces, for higher-order elements the difference in the number of DOFs can be quite significant. For instance, when comparing the finest discretization with $p=3$, the numbers of DOFs differ by a factor of approximately 2.6 (292,539 vs. 110,715 DOFs).

\subsection{Modal analysis -- cube}
We continue by performing a modal analysis in order to validate the proposed approach for dynamic problems. The computational domain is a cube of width 8
\begin{equation}\nonumber
    \Upomega=\left\{\left. (x, y, z) \in \mathbb{R}^3\,\right|\,0\leq x\leq 8,\,0\leq y\leq 8,\, 0\leq z\leq 8\right\}
\end{equation}
The material parameters are chosen as\\

\begin{tabular}{lc}
    \text{Young's modulus:} & $E =1$    \\
    \text{Poisson's ratio:} & $\nu =0$  \\
    \text{Mass density:}    & $\rho =1$ \\
\end{tabular}\\

\begin{figure}[t!]
    \centering
    \subfloat[Exemplary mesh, $h=4$, $p=3$, Mode~7 ($1^\text{st}$ nonzero mode). \label{fig:ModalConv_a}]{\includegraphics[width=0.28\textwidth]{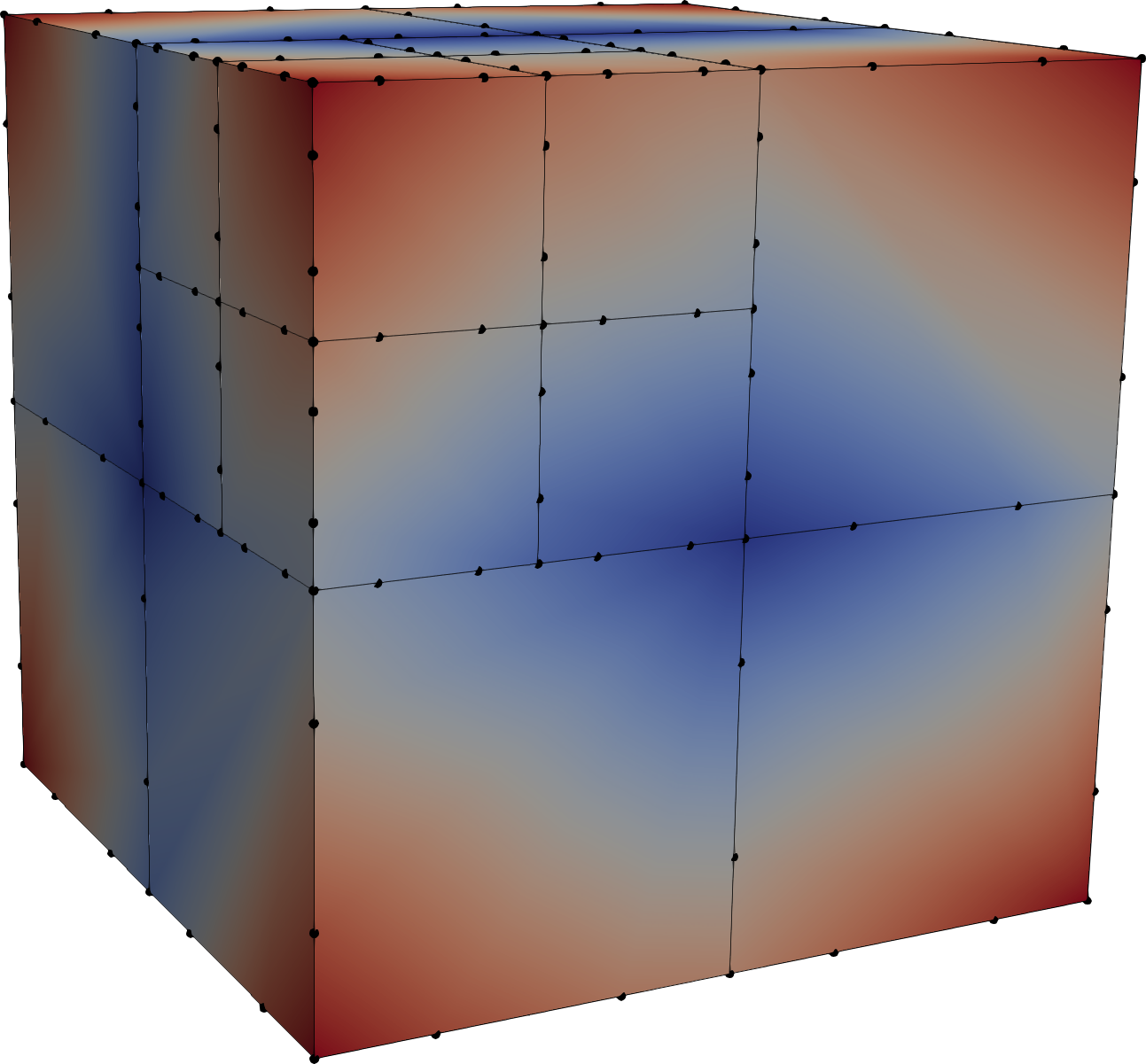}}\qquad \qquad
    \subfloat[Exemplary mesh, $h=2$, $p=1$, Mode~9 ($3^\text{rd}$ nonzero mode). \label{fig:ModalConv_b}]{\includegraphics[width=0.28\textwidth]{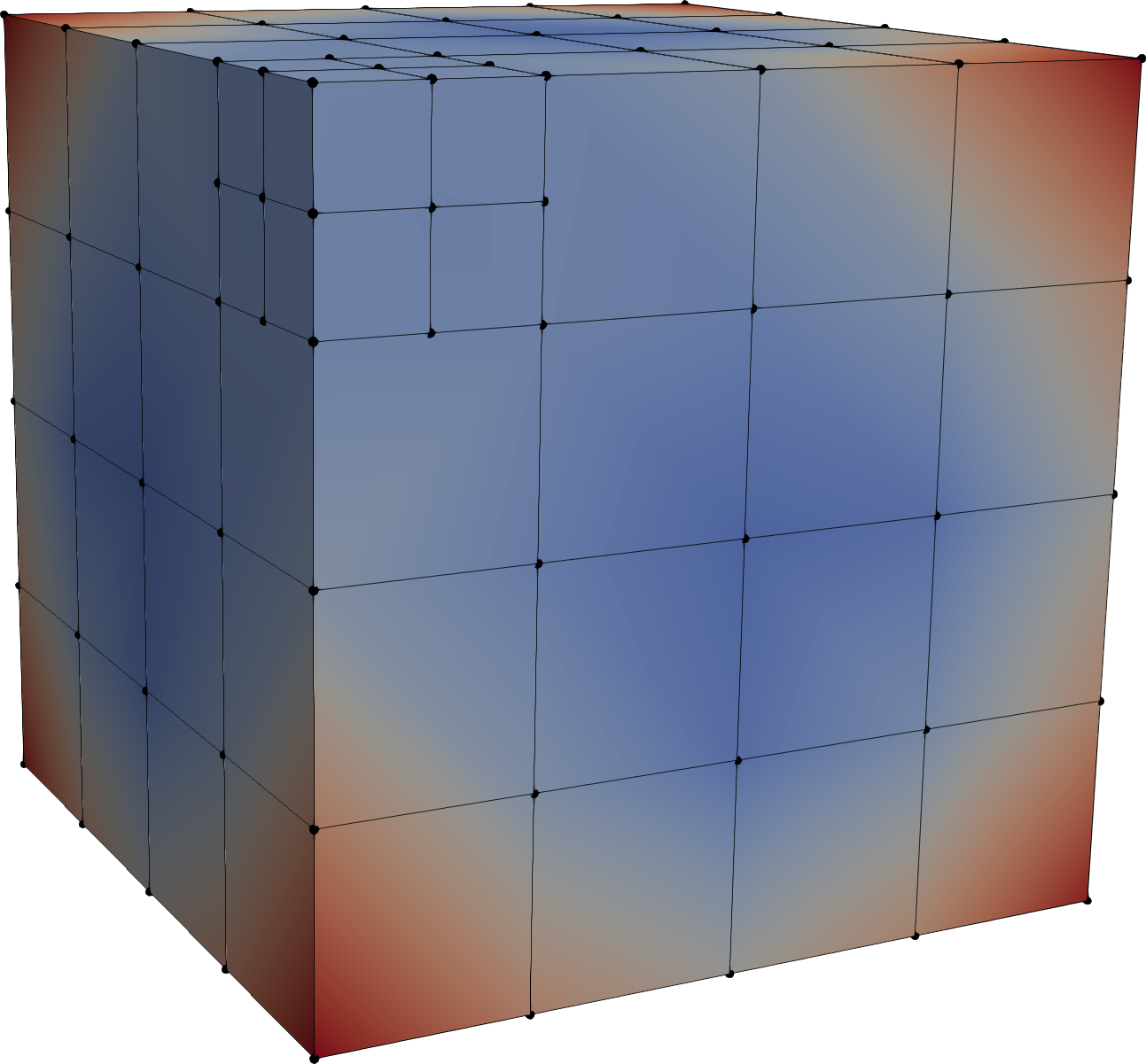}}\\
    \caption{Modal analysis of a cube -- Exemplary meshes and mode shapes.     \label{fig:ModalConv}}
\end{figure}

\begin{table}[t!]
    \centering
    \caption{Modal analysis of a cube -- Reference solution for the first ten nonzero eigenfrequencies.}
    \label{tab:modal_cube}
    \begin{tabular}{c|l}
        Mode No. & Eigenfrequency \\\hline
        7        & 0.063666938067 \\
        8        & 0.063666949380 \\
        9        & 0.108860021116 \\
        10       & 0.108860021166 \\
        11       & 0.108860027908 \\
        12       & 0.108860036839 \\
        13       & 0.108860080965 \\
        14       & 0.108861627176 \\
        15       & 0.117218751959 \\
        16       & 0.117218866414 \\
    \end{tabular}
\end{table}

\begin{figure}[b!]
    \centering
    \includegraphics[width=0.4\textwidth]{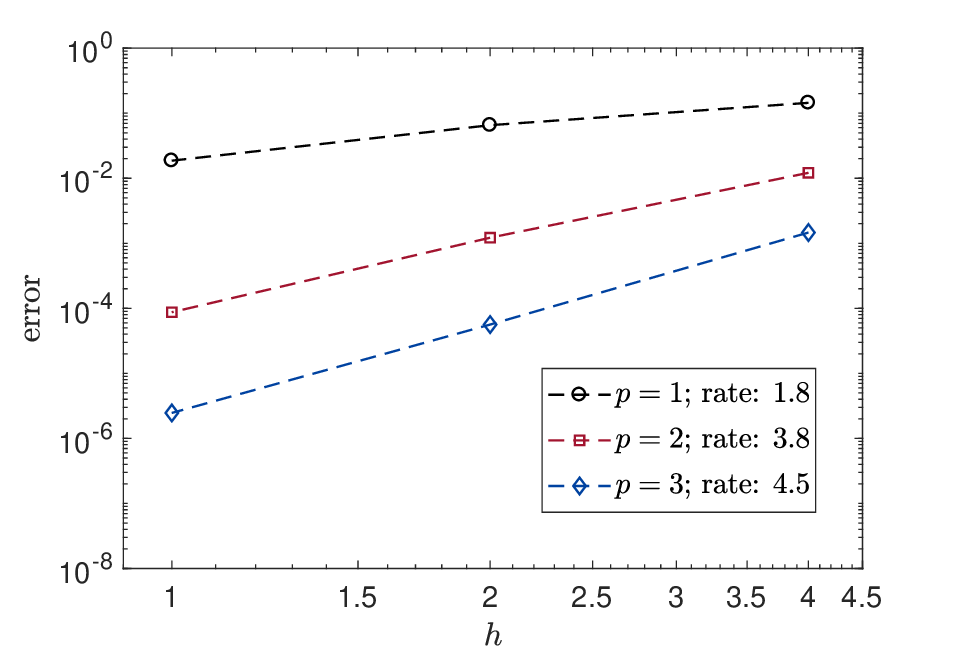}
    \caption{Modal analysis of a cube -- relative error of the eigenfrequencies.     \label{fig:ModalConv_error}}
\end{figure}

\noindent All surfaces are traction-free. The mesh is defined by $(8/h)^3$ cubic subdomains, where $h$ denotes the subdomains' side length.  To ensure that transition elements are present in the mesh, an SBFE-subdomain at one of the corners of the cube is divided into eight, see Fig.~\ref{fig:ModalConv} for examples of the mesh. A numerical reference solution has been computed by utilizing the conventional spectral element method on a uniform grid of $8^3$ elements of order 5, leading to 206,763 DOFs. The eigenfrequencies are listed in Table~\ref{tab:modal_cube}. The relative error of the first ten nonzero eigenfrequencies (modes 7-16) is plotted in Fig.~\ref{fig:ModalConv} with respect to the element size $h$ for different element orders $p$. When using $p=1,2,3$, the computed numerical convergence rate is obtained as 1.8, 3.8, and 4.5, respectively.

\subsection{Modal analysis -- Crane tower}
\begin{figure*}[p]
    \centering
    \subfloat{\includegraphics[width=0.45\textwidth]{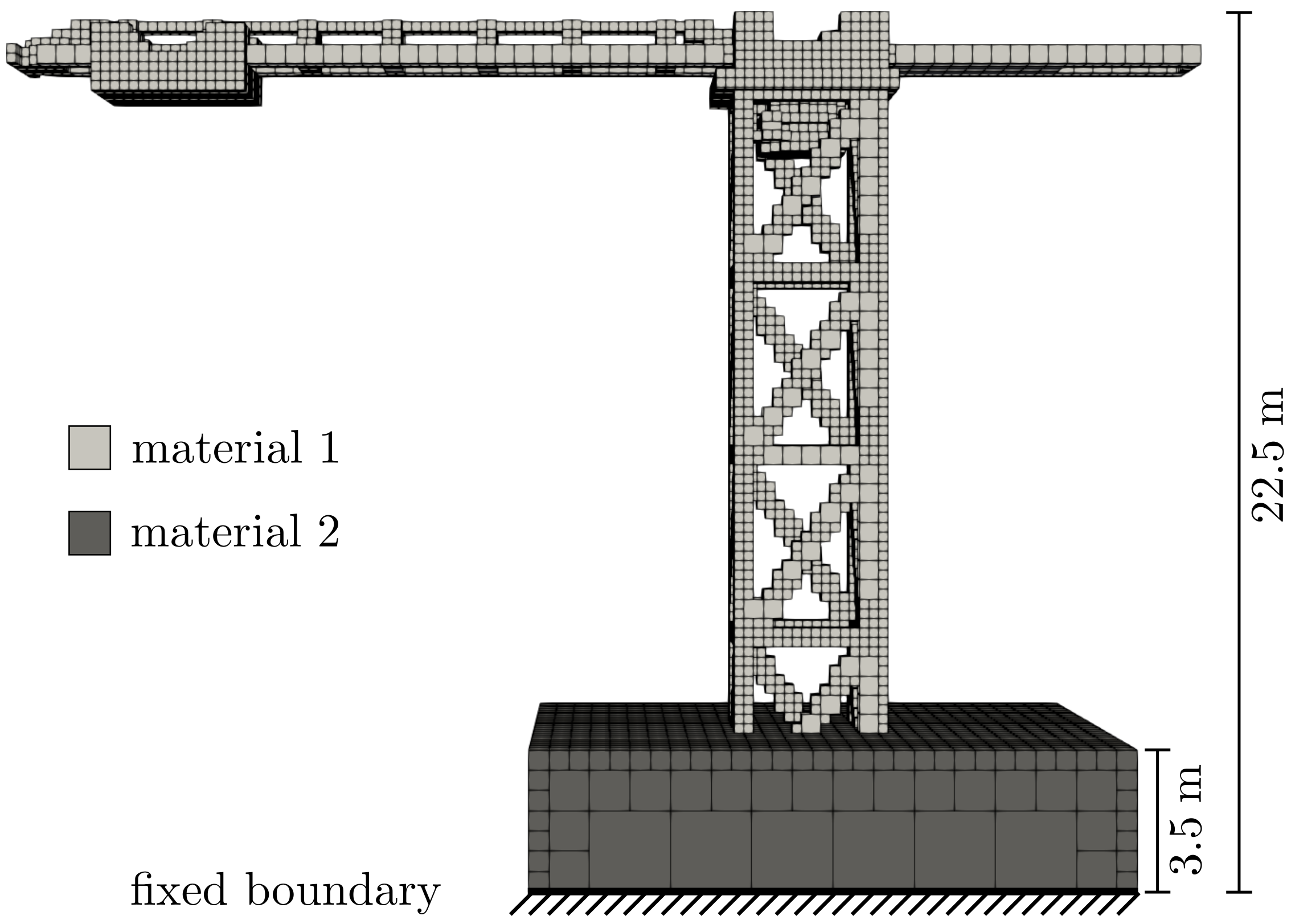}}\qquad
    \subfloat{\includegraphics[width=0.45\textwidth]{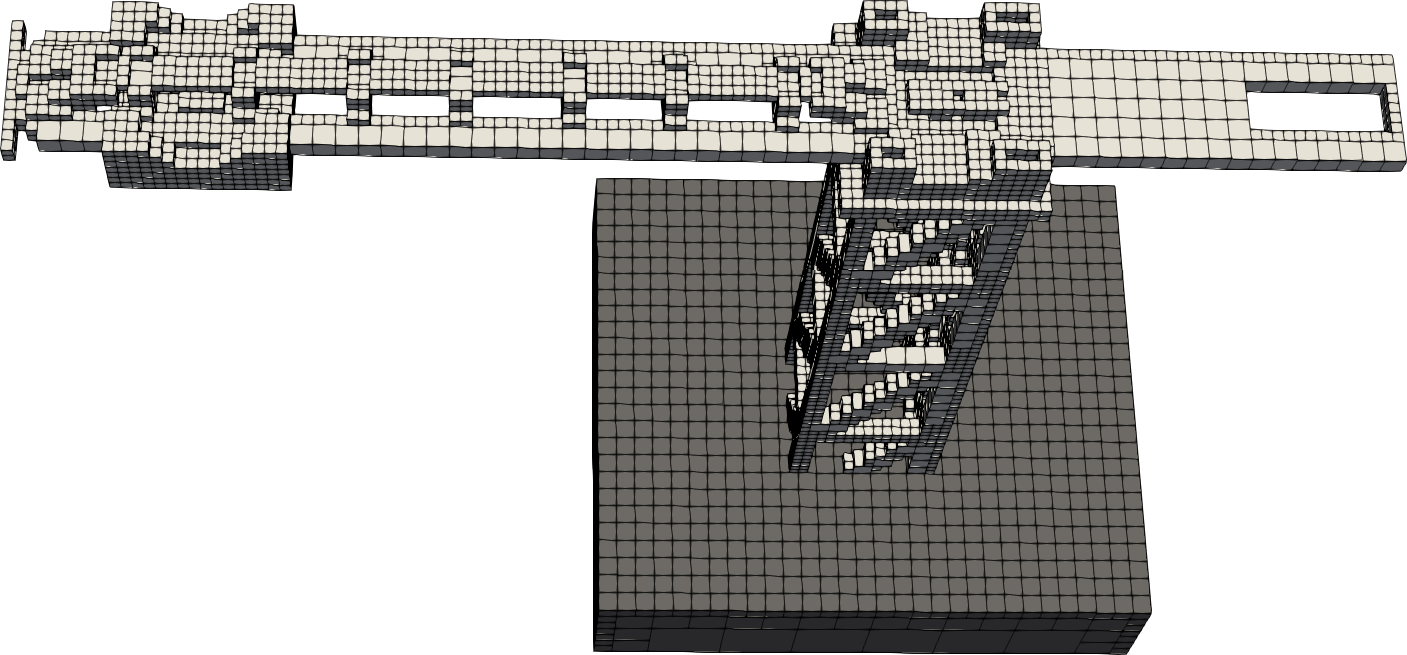}}
    \caption{Modal analysis of a crane tower -- octree decomposition.  \label{fig:crane_mesh}}
    \centering
    \subfloat[Mode 1, $f=1.1$ Hz]{\includegraphics[width=0.45\textwidth]{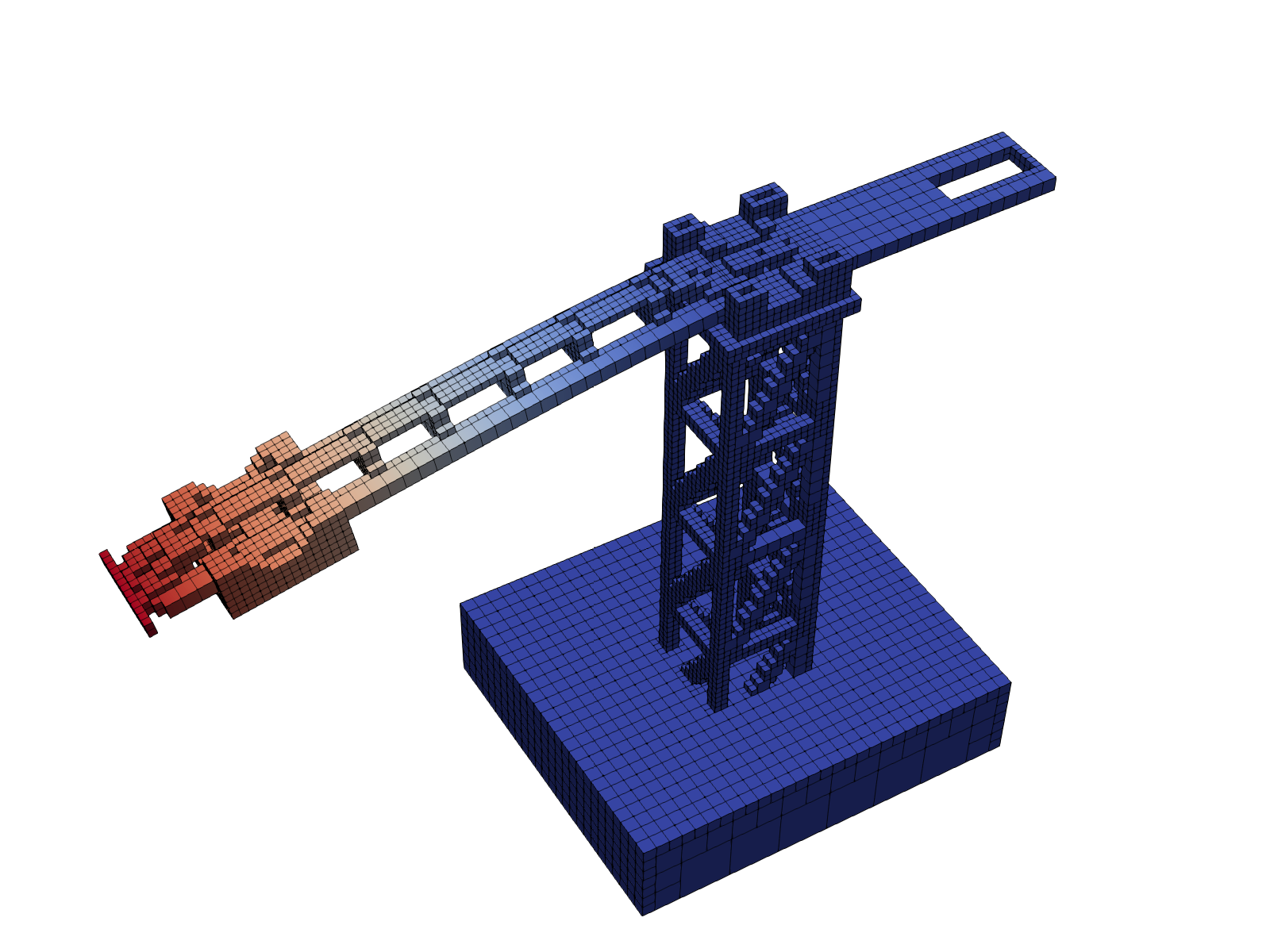}}
    \subfloat[Mode 2, $f=1.7$ Hz]{\includegraphics[width=0.45\textwidth]{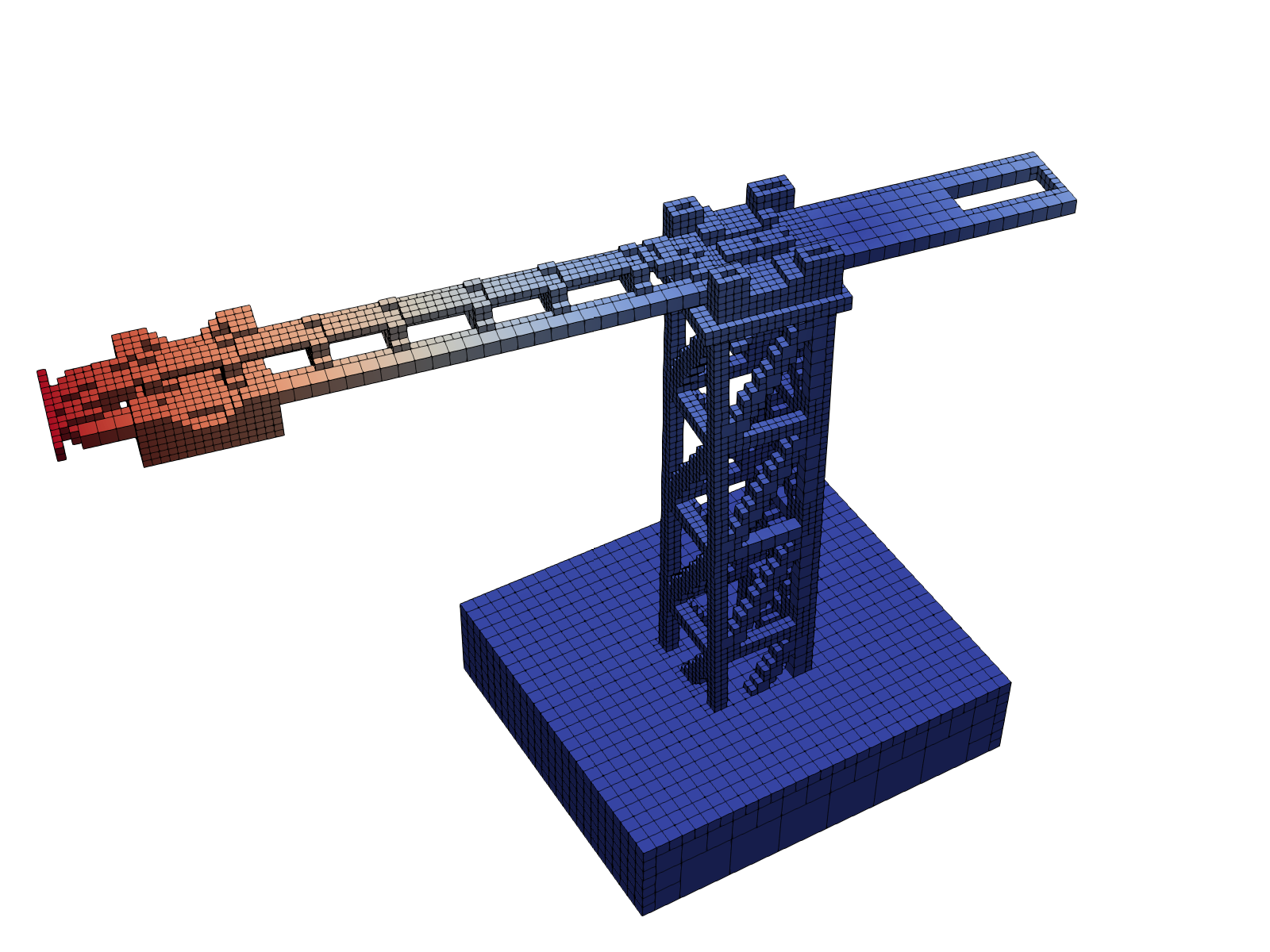}}\\
    \subfloat[Mode 3, $f=3.5$ Hz]{\includegraphics[width=0.45\textwidth]{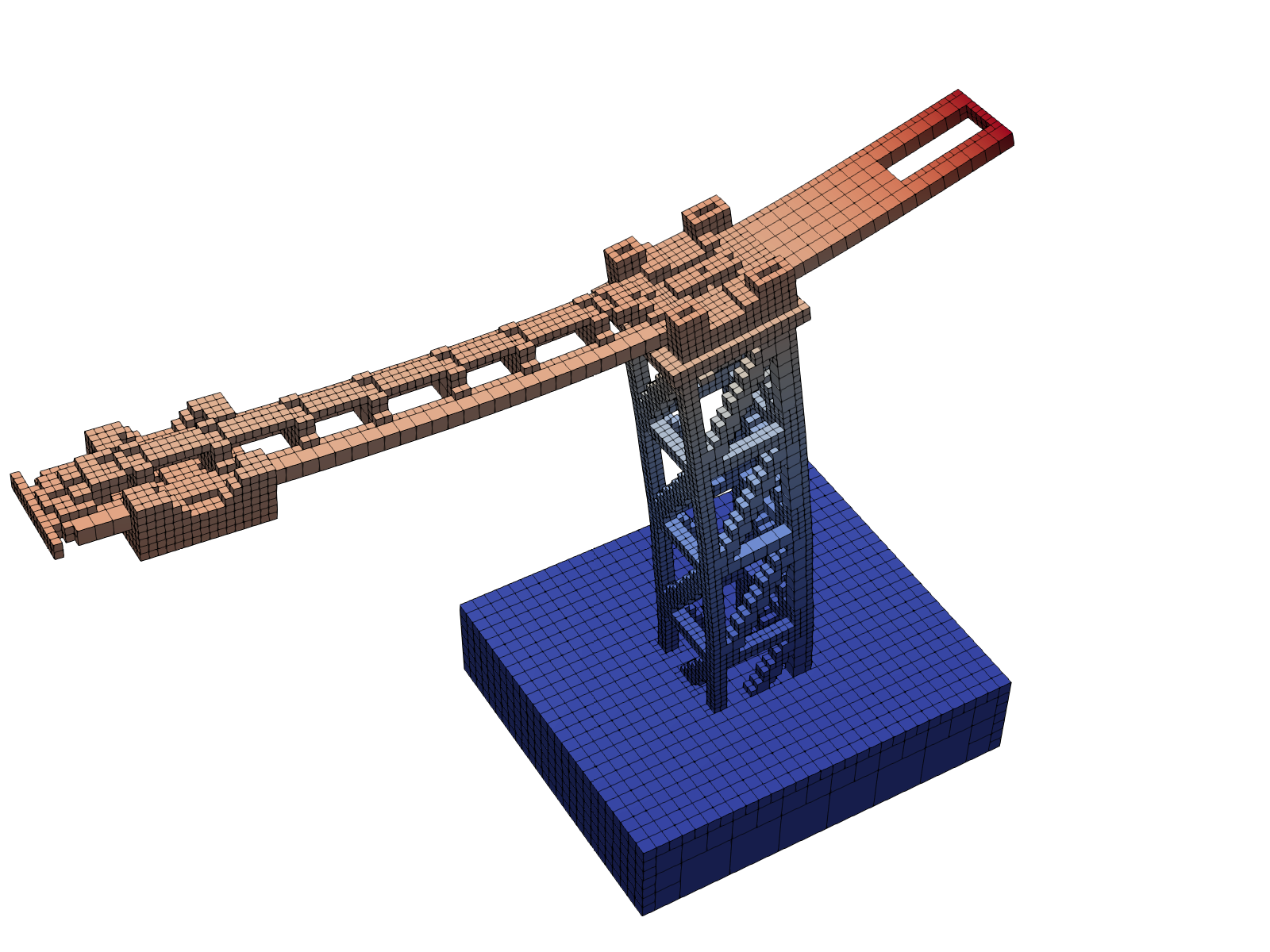}}
    \subfloat[Mode 4, $f=4.0$ Hz]{\includegraphics[width=0.45\textwidth]{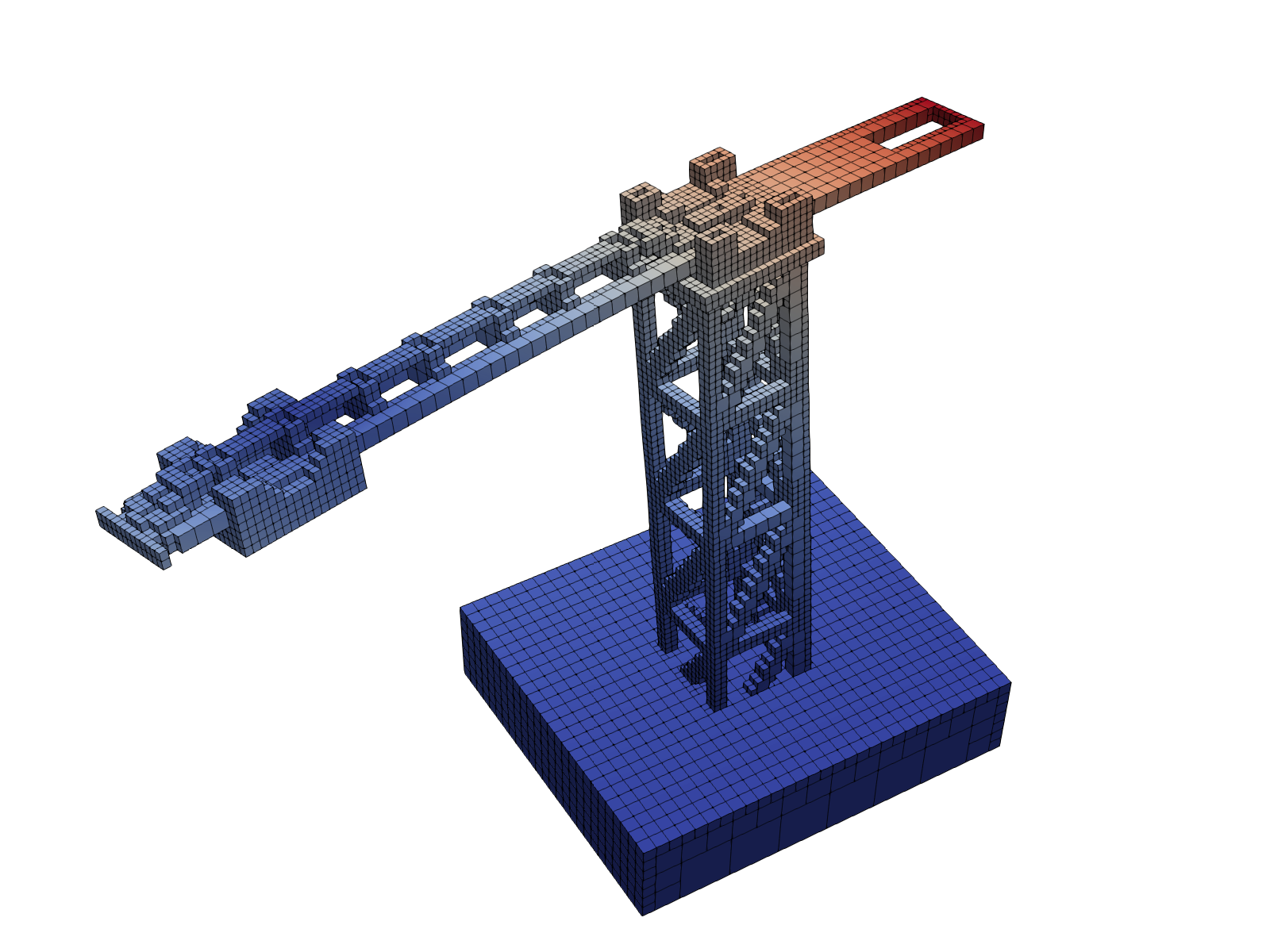}}
    \caption{Modal analysis of a crane tower -- Mode shapes and eigenfrequencies of the first four modes. \label{fig:crane_modes}}
\end{figure*}

To demonstrate the applicability of the proposed approach to more complex geometries, the model of a crane tower as depicted in Fig.~\ref{fig:crane_mesh} is analyzed. The model is based on a stereolithography (STL) file obtained from the online repository \textit{Thingiverse} \cite{MakerBotIndustries}.\footnote{The particular design \textit{crane tower} has been uploaded by Tung Nguyen (username \textit{tungnt}) and published under the terms of the GNU General Public License. The design can be retrieved from www.thingiverse.com/thing:2440385.}  From the STL-file, we generated a voxel-based model with the help of the online converter \textit{Voxelizer} \cite{Westerdiep2019}. We slightly modified the design -- mainly to remove unconstrained parts as well as to assign two different colors that correspond to different materials in the context of the image-based analysis. The modifications have been made using the software \textit{MagicaVoxel} \cite{Ephtracy2019}. The material parameters of the tower and the base (as indicated by the different colors) are assumed as \\

\begin{tabular}{lcc}
                                  & material 1   & material 2   \\
    \text{Young's modulus E: }    & 70 GPa       & 20 GPa       \\
    \text{Poisson's ratio $\nu$:} & 0.35         & 0.3          \\
    \text{mass density $\rho$:}   & 2.7 kg/m$^3$ & 2.4 kg/m$^3$ \\
\end{tabular}\\[0.5\baselineskip]

\noindent The total height of the model is 22.5\,m. We apply fixed boundary conditions at the bottom of the base, as indicated in Fig.~\ref{fig:crane_mesh}. The mesh has been created automatically by means of an octree decomposition of the voxel-based model, leading to 7158 subdomains. As can be inferred from Fig.~\ref{fig:crane_mesh}, the decomposition in this case leads to subdomains of four different sizes and therefore, the side length of the largest element is 8 times larger than that of the smallest elements. A rapid and consistent transition from large to small element sizes is easily achieved by the octree meshing algorithm. Figure~\ref{fig:crane_modes} shows the mode shapes of the first four modes and lists the corresponding eigenfrequencies. We found that to compute these modes, it suffices to employ quadratic elements on the surfaces of each subdomain, which results in a total of 615,960 DOFs.

\begin{figure*}[htb]
    \centering
    \subfloat[Exemplary mesh. \label{fig:castle_mesh}]{\includegraphics[width=0.45\textwidth]{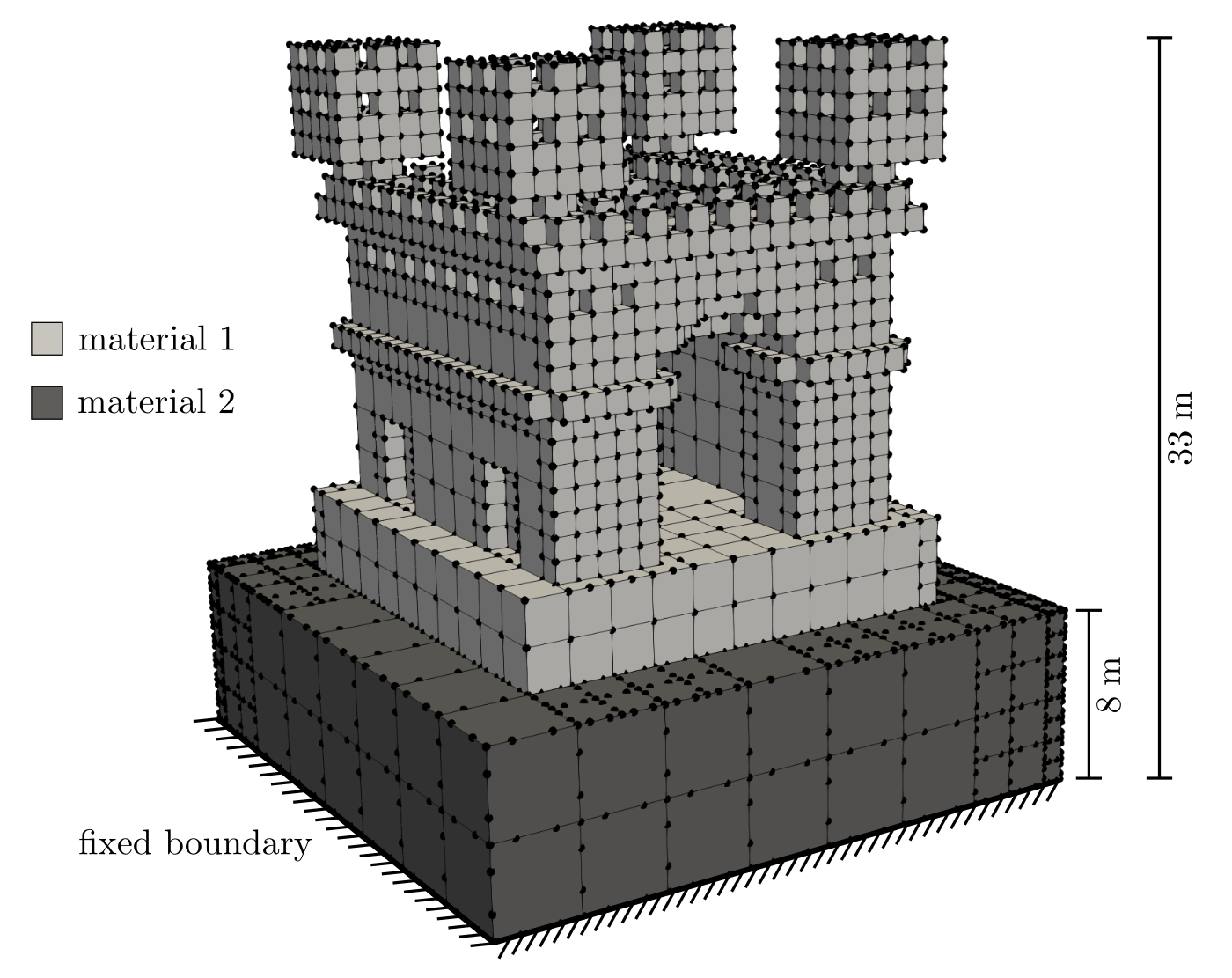}}
    \subfloat[Displacement field due to self-weight. \label{fig:castle_disp}]{\includegraphics[width=0.45\textwidth]{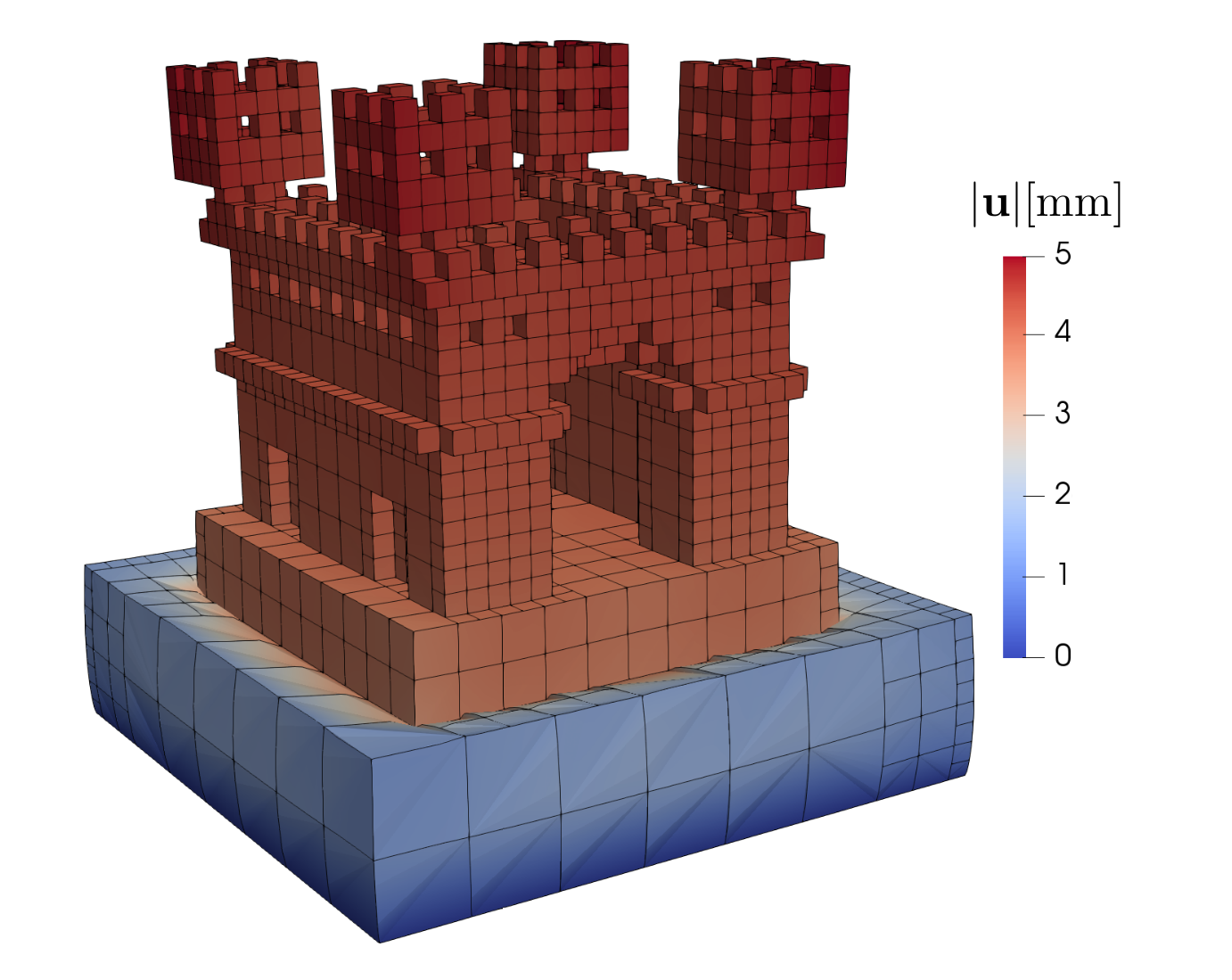}}
    \caption{Structure under self-weight loading -- (a) mesh and (b) displacement solution. \label{fig:castle}}
\end{figure*}

\subsection{Structure under self-weight loading}
As a final example, we study the behavior of the structure depicted in Fig.~\ref{fig:castle} under the influence of self-weight loading. The model of a castle is based on a sample provided in the software \textit{MagicaVoxel} \cite{Ephtracy2019} which we slightly modified and placed on a base of a different material. We scaled the model such that the total height of the structure is 33\,m. The base is a homogeneous layer with a thickness of 8\,m.
The material parameters of the castle (material 1) and foundation (material 2) are chosen as \\

\begin{tabular}{lcc}
                                  & material 1   & material 2 \\
    \text{Young's modulus E: }    & 10 GPa       & 0.5 GPa    \\
    \text{Poisson's ratio $\nu$:} & 0.3          & 0.2        \\
    \text{mass density $\rho$:}   & 2.4 kg/m$^3$ & 2 kg/m$^3$ \\
\end{tabular}\\

\noindent
The acceleration of gravity is assumed as $g=-9.81\,\text{m}/\text{s}^2$. As a proof of concept, this example represents a situation where it can be useful to apply a different element order to the two distinct materials and thus better capture the larger deformations in the softer material. This feature will be particularly interesting for more complicated problems, such as wave propagation through layered soils, which has been addressed previously in two dimensions \cite{Gravenkamp2017a}. In those cases, it can also be useful to adjust the element order not only based on the material properties but also on the size of the subdomain. In the current rather simple example, we found it sufficient to use elements of order 3 in the base material, while the comparably rigid material of the castle is discretized using linear elements only. The resulting mesh is depicted in Fig.~\ref{fig:castle_mesh}. The deformed geometry due to self-weight can be seen in Fig.~\ref{fig:castle_disp}, where the colors indicate absolute values (magnitude) of displacement.

\begin{table}[h]
    \centering
    \caption{Eigenfrequencies of the crane tower, computed using varying element order $p$.}
    \label{tab:castle_freqs}
    \begin{tabular}{c|rrr}
                                  & \multicolumn{3}{c}{frequency [Hz]}                                                         \\\hline
        \multicolumn{1}{c|}{Mode} & \multicolumn{1}{c}{$p=1$}          & \multicolumn{1}{c}{$p=2$} & \multicolumn{1}{c}{$p=3$} \\\hline
        1                         & 1.28                               & 1.12                      & 1.08                      \\
        2                         & 1.85                               & 1.74                      & 1.72                      \\
        3                         & 3.64                               & 3.49                      & 3.47                      \\
        4                         & 4.20                               & 4.02                      & 4.00                      \\
        5                         & 5.79                               & 5.02                      & 4.95                      \\
        6                         & 8.02                               & 7.24                      & 7.10                      \\
        7                         & 9.87                               & 8.85                      & 8.60                      \\
        8                         & 10.55                              & 10.12                     & 10.08                     \\
        9                         & 13.48                              & 12.81                     & 12.68                     \\
        10                        & 18.66                              & 18.43                     & 18.41                     \\
    \end{tabular}
\end{table}
%
%
\section{Conclusion}
The combination of the three-dimensional SBFEM with transition elements allows a consistent discretization of an octree decomposition. Each surface of each subdomain is discretized by quadrilateral elements. Hence, we avoid any unnecessary subdivision into smaller surface elements as had been done in previous publications. The numerical results demonstrate that the computational models created using the proposed technique pass the linear as well as higher order patch tests. Compared to the previous meshing paradigm (involving triangulation), there is no loss of accuracy due to the transition elements, while the number of surface elements and degrees of freedom is reduced. We also demonstrated that the proposed approach allows coupling elements of different interpolation orders straightforwardly. Thus, we are now able to implement local $p$-refinement, which could previously only be exploited in two-dimensional SBFEM models.
\FloatBarrier

\bibliographystyle{spmpsci}      
\bibliography{pnh_paper,paper_pnh_sbfem}

\begin{thebibliography}{10}
\providecommand{\url}[1]{{#1}}
\providecommand{\urlprefix}{URL }
\expandafter\ifx\csname urlstyle\endcsname\relax
  \providecommand{\doi}[1]{DOI~\discretionary{}{}{}#1}\else
  \providecommand{\doi}{DOI~\discretionary{}{}{}\begingroup
  \urlstyle{rm}\Url}\fi

\bibitem{Birk2012}
Birk, C., Prempramote, S., Song, C.: {An improved continued-fraction-based
  high-order transmitting boundary for time-domain analyses in unbounded
  domains}.
\newblock International Journal for Numerical Methods in Engineering
  \textbf{89}, 269--298 (2012)

\bibitem{Birk2009}
Birk, C., Song, C.: {A continued-fraction approach for transient diffusion in
  unbounded medium}.
\newblock Computer Methods in Applied Mechanics and Engineering \textbf{198},
  2576--2590 (2009)

\bibitem{Birkhoff1974}
Birkhoff, G., Cavendish, J.C., Gordon, W.J.: {Multivariant approximation by
  locally blended univariate interpolants}.
\newblock Proceedings of the National Academy of Sciences of the United States
  of America \textbf{71}(9), 3423--3425 (1974)

\bibitem{Bishop2014}
Bishop, J.E.: {A displacement-based finite element formulation for general
  polyhedra using harmonic shape functions}.
\newblock International Journal for Numerical Methods in Engineering
  \textbf{97}, 1--31 (2014)

\bibitem{PhDBroeker2001}
Br\"oker, H.: {I}ntegration von geometrischer {M}odellierung und {B}erechnung
  nach der \emph{p}-{V}ersion der {FEM}.
\newblock Berichte aus dem Bauwesen. Shaker Verlag (2001)

\bibitem{ArticleCavendish1975}
Cavendish, J.C.: {Local Mesh Refinement Using Rectangular Blended Finite
  Elements}.
\newblock Journal of Computational Physics \textbf{19}, 211--228 (1975)

\bibitem{Chen2014b}
Chen, X., Birk, C., Song, C.: {Numerical modelling of wave propagation in
  anisotropic soil using a displacement unit-impulse-response-based formulation
  of the scaled boundary finite element method}.
\newblock Soil Dynamics and Earthquake Engineering \textbf{65}, 243--255 (2014)

\bibitem{Chiong2014}
Chiong, I., Ooi, E.T., Song, C., Tin-Loi, F.: {Scaled boundary polygons with
  application to fracture analysis of functionally graded materials}.
\newblock International Journal for Numerical Methods in Engineering
  \textbf{98}, 562--589 (2014)

\bibitem{du2006recent}
Du, Q., Wang, D.: Recent progress in robust and quality {D}elaunay mesh
  generation.
\newblock Journal of Computational and Applied Mathematics \textbf{195}(1),
  8--23 (2006)

\bibitem{PhDDuczek2014}
Duczek, S.: Higher Order Finite Elements and the Fictitious Domain Concept for
  Wave Propagation Analysis.
\newblock VDI Fortschritt-Berichte Reihe 20 Nr. 458 (2014)

\bibitem{ArticleDuczek2019a}
Duczek, S., Gravenkamp, H.: Critical assessment of different mass lumping
  schemes for higher order serendipity finite elements.
\newblock Computer Methods in Applied Mechanics and Engineering \textbf{350},
  836--897 (2019)

\bibitem{ArticleDuczek2019b}
Duczek, S., Gravenkamp, H.: Mass lumping techniques in the spectral element
  method: On the equivalence of the row-sum, nodal quadrature, and diagonal
  scaling methods.
\newblock Computer Methods in Applied Mechanics and Engineering \textbf{353},
  516--569 (2019)

\bibitem{Duczek2019b}
Duczek, S., Saputra, A.A., Gravenkamp, H.: {High order transition elements: The
  xNy-element concept -- Part I: Statics}.
\newblock ArXiv e-prints \textbf{1909.04899}, 1--46 (2019)

\bibitem{Ephtracy2019}
Ephtracy: {MagicaVoxel 0.99.3a} (2019).
\newblock \urlprefix\url{https://ephtracy.github.io/}

\bibitem{Gabbert1999}
Gabbert, U., Graeff-Weinberg, K.: {Adaptive local-global analysis by pNh
  transition elements}.
\newblock Technische Mechanik \textbf{19}(2), 115--126 (1999)

\bibitem{ArticleGordon1971}
Gordon, W.J.: Blending-function methods of bivariate and multivariate
  interpolation and approximation.
\newblock SIAM Journal on Numerical Analysis \textbf{8}, 158--177 (1971)

\bibitem{ArtcileGordon1973a}
Gordon, W.J., Hall, C.A.: Construction of curvilinear co-ordinate systems and
  applications to mesh generation.
\newblock International Journal for Numerical Methods in Engineering
  \textbf{7}, 461--477 (1973)

\bibitem{Gordon1973a}
Gordon, W.J., Hall, C.A.: {Transfinite Element Methods: Blending-Function
  Interpolation over Arbitrary Curved Element Domains}.
\newblock Numerische Mathematik \textbf{21}, 109--129 (1973)

\bibitem{ArtcileGordon1973b}
Gordon, W.J., Hall, C.A.: Transfinite element methods: Blending-function
  interpolation over arbitrary curved element domains.
\newblock Numerische Mathematik \textbf{21}, 109--129 (1973)

\bibitem{BookGordon1982}
Gordon, W.J., Thiel, L.C.: Transfinite Mappings and their Application to Grid
  Generation.
\newblock Elsevier Science Publishing (1982)

\bibitem{Graeff-Weinberg1996a}
Graeff-Weinberg, K., Berger, H.: {Verbesserte FE-Diskretisierung bei
  Kontaktaufgaben}.
\newblock Technische Mechanik \textbf{16}(3), 250--270 (1996)

\bibitem{Gravenkamp2018}
Gravenkamp, H.: {Efficient simulation of elastic guided waves interacting with
  notches, adhesive joints, delaminations and inclined edges in plate
  structures}.
\newblock Ultrasonics \textbf{82}, 101--113 (2018)

\bibitem{Gravenkamp2014a}
Gravenkamp, H., Bause, F., Song, C.: {On the computation of dispersion curves
  for axisymmetric elastic waveguides using the scaled boundary finite element
  method}.
\newblock Computers {\&} Structures \textbf{131}, 46--55 (2014)

\bibitem{Gravenkamp2014f}
Gravenkamp, H., Birk, C., Song, C.: {Simulation of elastic guided waves
  interacting with defects in arbitrarily long structures using the scaled
  boundary finite element method}.
\newblock Journal of Computational Physics \textbf{295}, 438--455 (2015)

\bibitem{Gravenkamp2017c}
Gravenkamp, H., Duczek, S.: {Automatic image-based analyses using a coupled
  quadtree-SBFEM/SCM approach}.
\newblock Computational Mechanics \textbf{60}, 559--584 (2017)

\bibitem{Gravenkamp2018a}
Gravenkamp, H., Natarajan, S.: {Scaled boundary polygons for linear
  elastodynamics}.
\newblock Computer Methods in Applied Mechanics and Engineering \textbf{333},
  238--256 (2018)

\bibitem{Gravenkamp2019}
Gravenkamp, H., Saputra, A.A., Duczek, S.: {High-order shape functions in the
  scaled boundary finite element method revisited}.
\newblock Archives of Computational Methods in Engineering \textbf{under
  review}, 1--21 (2019)

\bibitem{Gravenkamp2017a}
Gravenkamp, H., Saputra, A.A., Song, C., Birk, C.: {Efficient wave propagation
  simulation on quadtree meshes using SBFEM with reduced modal basis}.
\newblock International Journal for Numerical Methods in Engineering
  \textbf{110}, 1119--1141 (2017)

\bibitem{MakerBotIndustries}
Industries, M.: {Thingiverse}.
\newblock \urlprefix\url{https://www.thingiverse.com/}

\bibitem{BookKarniadakis2005}
Karniadakis, G.E., Sherwin, S.J.: Spectral/\emph{hp} Element Methods for
  Computational Fluid Dynamics.
\newblock Oxford Science Publications (2005)

\bibitem{Kausel1994b}
Kausel, E.: {Thin-layer method: Formulation in the time domain}.
\newblock International Journal for Numerical Methods in Engineering
  \textbf{37}, 927--941 (1994)

\bibitem{Kausel2019}
Kausel, E., Gravenkamp, H.: {On the numerical solution of matrix Bessel
  equations}.
\newblock ZAMM Zeitschrift fur Angewandte Mathematik und Mechanik
  \textbf{99}(8), e201800288 (2019)

\bibitem{Kausel1981b}
Kausel, E., Ro{\"{e}}sset, J.M.: {Stiffness matrices for layered soils}.
\newblock Bulletin of the Seismological Society of America \textbf{71}(6),
  1743--1761 (1981)

\bibitem{keyak1990automated}
Keyak, J., Meagher, J., Skinner, H., Mote, C.: Automated three-dimensional
  finite element modelling of bone: a new method.
\newblock Journal of Biomedical Engineering \textbf{12}(5), 389--397 (1990)

\bibitem{ArticleKiralyfalvi1997}
Kir\'{a}lyfalvi, G., Szab\'{o}, B.: Quasi-regional mapping for the
  \emph{p}-version of the finite element method.
\newblock Finite Elements in Analysis and Design \textbf{27}, 85--97 (1997)

\bibitem{Krome2017}
Krome, F., Gravenkamp, H.: {A semi-analytical curved element for linear
  elasticity based on the scaled boundary finite element method}.
\newblock International Journal for Numerical Methods in Engineering
  \textbf{109}, 790--808 (2017)

\bibitem{Krome2017a}
Krome, F., Gravenkamp, H., Birk, C.: {Prismatic semi-analytical elements for
  the simulation of linear elastic problems in structures with piecewise
  uniform cross section}.
\newblock Computers {\&} Structures \textbf{192}, 83--95 (2017)

\bibitem{Liu2019a}
Liu, L., Zhang, J., Song, C., Birk, C., Gao, W.: {An automatic approach for the
  acoustic analysis of three-dimensional bounded and unbounded domains by
  scaled boundary finite element method}.
\newblock International Journal of Mechanical Sciences \textbf{151}, 563--581
  (2019)

\bibitem{Liu2017}
Liu, Y., Saputra, A.A., Wang, J., Tin-Loi, F., Song, C.: {Automatic polyhedral
  mesh generation and scaled boundary finite element analysis of STL models}.
\newblock Computer Methods in Applied Mechanics and Engineering \textbf{313},
  106--132 (2017)

\bibitem{lohner1988generation}
L{\"o}hner, R., Parikh, P.: Generation of three-dimensional unstructured grids
  by the advancing-front method.
\newblock International Journal for Numerical Methods in Fluids \textbf{8}(10),
  1135--1149 (1988)

\bibitem{lorensen1987marching}
Lorensen, W.E., Cline, H.E.: Marching cubes: A high resolution 3d surface
  construction algorithm.
\newblock In: ACM SIGGRAPH Computer Graphics, vol.~21, pp. 163--169. ACM (1987)

\bibitem{Man2012}
Man, H., Song, C., Gao, W., Tin-Loi, F.: {A unified 3D-based technique for
  plate bending analysis using scaled boundary finite element method}.
\newblock International Journal for Numerical Methods in Engineering
  \textbf{91}, 491--515 (2012)

\bibitem{Man2014}
Man, H., Song, C., Natarajan, S., Ooi, E.T., Birk, C., {Tat Ooi}, E., Birk, C.:
  {Towards automatic stress analysis using Scaled Boundary Finite Element
  Method with quadtree mesh of high-order elements}.
\newblock ArXiv e-prints p. math.NA/1402.5186 (2014)

\bibitem{Ooi2015b}
Ooi, E.T., Man, H., Natarajan, S., Song, C.: {Adaptation of quadtree meshes in
  the scaled boundary finite element method for crack propagation modelling}.
\newblock Engineering Fracture Mechanics \textbf{144}, 101--117 (2015)

\bibitem{Ooi2013}
Ooi, E.T., Shi, M., Song, C., Tin-Loi, F., Yang, Z.: {Dynamic crack propagation
  simulation with scaled boundary polygon elements and automatic remeshing
  technique}.
\newblock Engineering Fracture Mechanics \textbf{106}(2012), 1--21 (2013)

\bibitem{Ooi2014}
Ooi, E.T., Song, C., Tin-Loi, F.: {A scaled boundary polygon formulation for
  elasto-plastic analyses}.
\newblock Computer Methods in Applied Mechanics and Engineering \textbf{268},
  905--937 (2014)

\bibitem{Ooi2012a}
Ooi, E.T., Song, C., Tin-Loi, F., Yang, Z.: {Automatic modelling of cohesive
  crack propagation in concrete using polygon scaled boundary finite elements}.
\newblock Engineering Fracture Mechanics \textbf{93}, 13--33 (2012)

\bibitem{BookPozrikidis2014}
Pozrikidis, C.: Introduction to Finite and Spectral Element Methods using
  MATLAB, 2 edn.
\newblock Chapman and Hall/CRC (2014)

\bibitem{ArticleProvatidis2006}
Provatidis, C.G.: Coons-patch macroelements in two-dimensional parabolic
  problems.
\newblock Applied Mathematical Modelling \textbf{30}(4), 319--351 (2006)

\bibitem{ArticleProvatidis2011}
Provatidis, C.G.: Two-dimensional elastostatic analysis using {C}oons-{G}ordon
  interpolation.
\newblock Meccanica \textbf{47}(4), 951--967 (2011)

\bibitem{BookProvatidis2019}
Provatidis, C.G.: Precursors of Isogeometric Analysis.
\newblock Springer International Publishing (2019)

\bibitem{Saputra2015}
Saputra, A.A., Birk, C., Song, C.: {Computation of three-dimensional fracture
  parameters at interface cracks and notches by the scaled boundary finite
  element method}.
\newblock Engineering Fracture Mechanics \textbf{148}, 213--242 (2015)

\bibitem{Saputra2017}
Saputra, A.A., Talebi, H., Tran, D., Birk, C., Song, C.: {Automatic image-based
  stress analysis by the scaled boundary finite element method}.
\newblock International Journal for Numerical Methods in Engineering
  \textbf{109}, 697--738 (2017)

\bibitem{Song2004}
Song, C.: {A matrix function solution for the scaled boundary finite-element
  equation in statics}.
\newblock Computer Methods in Applied Mechanics and Engineering \textbf{193},
  2325--2356 (2004)

\bibitem{Song2009}
Song, C.: {The scaled boundary finite element method in structural dynamics}.
\newblock International Journal for Numerical Methods in Engineering
  \textbf{77}, 1139--1171 (2009)

\bibitem{Song2018}
Song, C.: {The scaled boundary finite element method: introduction to theory
  and implementation}.
\newblock Wiley (2018)

\bibitem{Song1997}
Song, C., Wolf, J.P.: {The scaled boundary finite-element method - alias
  consistent infinitesimal finite-element cell method - for elastodynamics}.
\newblock Computer Methods in Applied Mechanics and Engineering \textbf{147},
  329--355 (1997)

\bibitem{Song2000a}
Song, C., Wolf, J.P.: {The scaled boundary finite-element method - a primer:
  solution procedures}.
\newblock Computers {\&} Structures \textbf{78}, 211--225 (2000)

\bibitem{BookSzabo1991}
Szab\'{o}, B., Babu\v{s}ka, I.: Finite Element Analysis.
\newblock John Wiley and Sons (1991)

\bibitem{Timoshenko1951}
Timoshenko: {Theory of elasticity} (1951)

\bibitem{PhDWeinberg1996}
Weinberg, K.: {E}in {F}inite-{E}lemente-{K}onzept zur lokalen {N}etzverdichtung
  und seine {A}nwendung auf {K}oppel- und {K}ontaktprobleme.
\newblock Ph.D. thesis, Otto von Guericke University Magdeburg (1996)

\bibitem{ArticleWeinberg2002}
Weinberg, K., Gabbert, U.: An adaptive \emph{p{N}h}-technique for global-local
  finite element analysis.
\newblock Engineering Computations \textbf{19}, 485--500 (2002)

\bibitem{Weinberg2002}
Weinberg, K., Gabbert, U.: {An adaptive pNh-technique for global-local finite
  element analysis}.
\newblock Engineering Computations \textbf{19}(5), 485--500 (2002)

\bibitem{Westerdiep2019}
Westerdiep, A.: {Online Voxelizer} (2019).
\newblock \urlprefix\url{http://drububu.com/miscellaneous/voxelizer/}

\bibitem{Wolf1994}
Wolf, J.P., Song, C.: {Dynamic-stiffness matrix in time domain of unbounded
  medium by infinitesimal finite element cell method}.
\newblock Earthquake Engineering {\&} Structural Dynamics \textbf{23},
  1181--1198 (1994)

\bibitem{Wolf1994a}
Wolf, J.P., Song, C.: {Dynamic-stiffness matrix of unbounded soil by
  finite-element multi-cell cloning}.
\newblock Earthquake Engineering {\&} Structural Dynamics \textbf{23}, 233--250
  (1994)

\bibitem{Wolf2000a}
Wolf, J.P., Song, C.: {The scaled boundary finite-element method - a primer:
  derivations}.
\newblock Computers {\&} Structures \textbf{78}, 191--210 (2000)

\bibitem{young2008efficient}
Young, P., Beresford-West, T., Coward, S., Notarberardino, B., Walker, B.,
  Abdul-Aziz, A.: An efficient approach to converting three-dimensional image
  data into highly accurate computational models.
\newblock Philosophical Transactions of the Royal Society A: Mathematical,
  Physical and Engineering Sciences \textbf{366}(1878), 3155--3173 (2008)

\end{thebibliography}

\end{document}